\documentclass[10pt]{article}
\usepackage{amssymb}
\newtheorem{theorem}{Theorem}[section]
\newtheorem{corollary}{Corollary}[section]
\newtheorem{lemma}{Lemma}[section]

\newtheorem{conjecture}{Conjecture}[section]

\newtheorem{remark}{Remark}[section]
\newtheorem{definition}{Definition}[section]
\newtheorem{proposition}{Proposition}[section]

\newcommand{\halmos}{\rule{1ex}{1.4ex}}
				
\newcommand{\bea}{\begin{eqnarray}}
\newcommand{\epf}{\hspace*{\fill}\mbox{$\halmos$}}
\newcommand{\eea}{\end{eqnarray}}
\newcommand{\nn}{\nonumber \\}
\newcommand{\be}{\begin {equation}}

\newcommand{\ee}{\end{equation}}

\title{{Fusion rings for degenerate minimal
models}}
\author{Antun Milas}
\date{}
\begin{document}
\bibliographystyle{alpha}
\maketitle
\begin{abstract}
We study fusion rings for degenerate minimal models ($p=q$ case) for
$N=0$ and $N=1$ (super)conformal algebras.
We consider a distinguished 
family of modules at the level $c=1$ and $c=\frac{3}{2}$ 
and show that the corresponding fusion rings are isomorphic
to the representation rings for ${\mathfrak sl}(2, {\bf C})$
and ${\mathfrak osp}(1|2)$ respectively.

\end{abstract}

\section{Introduction}
The Virasoro algebra and its minimal models are a good source of
interesting vertex operator algebras.
In \cite{W} the rationality of
the Virasoro vertex operator algebras $L(c_{p.q},0)$ was proved, where
$c_{p,q}=1-\frac{6(p-q)^2}{pq}$ and $(p,q)=1$, $p,q \geq 2$. 
This result is used for the construction of the corresponding
vertex tensor categories (cf. \cite{H1}).
A similar result is obtained for $N=1$ case in \cite{Ad} and \cite{HM}. 

In this paper we study a non--rational
vertex operator algebra $L(1,0)$ ($p=q$ case)
and the corresponding fusion ring
for {\em degenerate} minimal models, i.e., the case ``$p=q$'',
with central charge $c=1$. We also consider a $N=1$ vertex operator superalgebra
version based on $L(\frac{3}{2},0)$ (see below).
These cases are substantially different for many reasons (let us focus on
the case $L(1,0)$ since the same problem persists for $L(\frac{3}{2},0)$.
The vertex operator algebra $L(1,0)$ is not {\em rational} (cf. \cite{FZ})
but it has a distinguished family of irreducible modules (those that
are not irreducible Verma modules) ${\cal F}_1$, 
which consists of classes of irreducible modules 
isomorphic to $L(1,\frac{m^2}{4})$ for some $m \in {\mathbb N}$.
These modules have a quite simple embedding structure (\cite{KR},
\cite{FF2}).


We show that
the fusion ring for the family ${\cal F}_1$ is 
isomorphic to the representation ring 
${\cal R}ep ({\mathfrak sl}(2, {\mathbb C}))$, i.e.,  
we ``formally'' have 
$$L\left(1,\frac{n^2}{4}\right) \times
L\left(1,\frac{m^2}{4}\right)=$$
$$L\left(1,\frac{(n+m)^2}{4}\right) +
L\left(1,\frac{(n+m-2)^2}{4}\right)  \ldots
+L\left(1,\frac{(n-m)^2}{4}\right),$$
where $m,n \in {\mathbb N}$ and $n \geq m$.  

This result seems to be known--in some form-- for a while by physicists 
(also in \cite{FKRW} is stated as a part
of more general conjecture concerning fusion rings 
for $W({\mathfrak gl}_N)$ algebras--see also \cite{FM}).
The author of the current paper could not trace any 
proof in the language of vertex operator algebras.
Some computations are done in \cite{DG} but not complete.
But instead of trying to patch missing proofs, there are
two more important reasons for seeking
such a proof.
\begin{itemize}
\item 
So far, not many computations of the fusion coefficients 
has been known for non--rational vertex operator
algebras (here non--rational means {\em non--rational} in any
reasonable category). In particular we offer a proof that
uses universal construction (induced modules), therefore
it is very general. 

\item As noticed by H. Li in \cite{L1} and \cite{L2}, Frenkel--Zhu's formula
\cite{FZ} does not hold for non--rational vertex operator algebras.
The right formula was provided in \cite{L2} but it is a non--trivial
matter to use it for computational purposes in non--rational setting.

\end{itemize}
We believe that our method can be used for more complicated models--like
degenerate models associated to ${\mathcal W}$--algebras.

We have to stress that the fusion coefficients 
are simply derived from the space of intertwining
operators among irreducible modules.
In other words it is {\em not} true
that the only modules which ``fuse'' with 
$L\left(1,\frac{n^2}{4}\right)$ and $L\left(1,\frac{m^2}{4}\right)$
are completely reducible.
This fact makes impossible to
implement $P(z)$--tensor product construction
from \cite{HL1}--\cite{HL2}. The resolution
might be to construct (a new) tensor product
which takes only irreducible modules into account, but
this approach will assume a good knowledge of matrix
coefficients for product of intertwining operators.
A different approach would be working in the
larger family $\bar{\cal F}_1$, which consists
of all quotients of Verma modules $M(1,\frac{m^2}{4})$. 
The possible constructions will be discussed elsewhere. 

We also provide a different proof of the fusion formulas
by constructing all intertwining operators from the 
lattice vertex operator algebra $V_L$ and
its irreducible module $V_{L+1/2}$ (cf. \cite{DG}).

A super $N=1$ versions of the above result stems from
the $N=1$ Neveu--Schwarz Lie superalgebra at the level $\frac{3}{2}$.
Again, there are essentially two approaches: one which 
uses the lattice construction (extended with a
suitable fermionic Fock space) and the other which uses
the singular vectors and projection formulas.
For the future purposes we use the latter approach. 
We consider a set of equivalence classes of irreducible 
modules for the $N=1$ superconformal algebra (see Section 3.) with 
representatives $L(\frac{3}{2}, \frac{q^2}{2})$ where $q \in \mathbb{N}.$
We proved (see Theorem \ref{last} and Corollary \ref{last1}) 
that the corresponding fusion ring is isomorphic to the representation ring
for $\mathfrak{osp}(1|2)$, i.e., we
formally have:
$$ L\left(\frac{3}{2},\frac{r^2}{2}\right) \times
L\left(\frac{3}{2},\frac{q^2}{2}\right)=$$
$$L\left(\frac{3}{2},\frac{(r+q)^2}{2}\right) +
L\left(\frac{3}{2},\frac{(r+q-1)^2}{2}\right)  \ldots
+L\left(\frac{3}{2},\frac{(r-q)^2}{2}\right),$$
for every $r,q \in {\bf N}$, $r \geq q$,
where $\times$ stands for the fusion product (see the last Chapter).

In particular, as in the Virasoro algebra case, these fusion coefficients
are $0$ or $1$. However in \cite{HM} we showed that for $N=1$ case 
has some interesting features; for some vertex operator algebras $L(c,0)$, 
fusion coefficients might be $2$. In Proposition \ref{twons}
we construct a non--trivial example with $c=\frac{15}{2}-3\sqrt{5}$.

At the very end, we construct an example of a {\em logarithmic intertwining 
operator} (for the definition see \cite{M})  for the $N=1$ vertex operator 
superalgebra  $L(\frac{27}{2},0)$.

\noindent {\em n.b.} These results can be extended for a more general class
of vertex operator algebras $L(c,0)$ where $c \neq c_{p,q}$; 
because of simplicity we treat only the case $c=1$ and $c=\frac{3}{2}$.

{\bf Acknowledgment:} The author thank Prof. Haisheng Li and Prof. Yi-Zhi Huang
for useful comments. Thanks go to the referee for his/her valuable remarks.
This paper is a union of slightly modified preprints available at the arXive.

\section{Representations of the Virasoro algebra at the level $c=1$}

The representation theory for the Virasoro algebra has been studied intensively
in the last two decades (\cite{KR}, \cite{FF1}--\cite{FF3}). Kac's determinant
formula is the most important tool in the highest (or lowest) 
weight theory.
From the determinant formula it follows that
the lowest weight Verma module with the central
charge $c(t)=13-6t-6t^{-1}$ and the weight 
$$h_{p,q}(t)=\frac{1-p^2}{4}t^{-1}-\frac{1-pq}{2}+\frac{1-q^2}{4}t,$$
has a singular vector of the weight $h_{p.q}(t)+pq$, $t \in {\mathbb C}$.
We are interested in the case $t=1$, i.e $c=1$. 
It is easy to see that $M(1,h)$ is irreducible if and only if 
$h \neq \frac{m^2}{4}$ for some $m \in {\mathbb N}$. In the case
$h=\frac{m^2}{4}$ we have the following description:
\begin{proposition} \label{sequence}
The Verma module $M(1,\frac{m^2}{4})$ has a unique singular vector
of  weight $\frac{m^2}{4}+(m+1)$. This vector generates the
maximal submodule. In other words we
have the following exact sequence
\begin{equation}  \label{exact}
0 \rightarrow M\left(1,\frac{(m+2)^2}{4}\right) \rightarrow M\left(1,\frac{m^2}{4}\right) \rightarrow
L\left(1,\frac{m^2}{4}\right) \rightarrow 0. 
\end{equation}
\end{proposition}

Even though they do not exist in general, 
in the case 
$h_{1,q}(t)$, if $p=1$ there are explicit
formulas  at each level $c(t)$ (in particular  $t=1$).
When $c=1$ Benoit and S. Aubin's formula \cite{BSA0}
implies that
\begin{equation} \label{singular}
P_{{\rm sing}}v_{1,q}=\sum_{\stackrel{ I=\{i_1,...,i_n \}
}{|I|=q}}c_{q}(i_1,...,i_n)L(-i_1)\ldots L(-i_n)v_{1.q}
\end{equation}
is a singular vector for $M(1,h_{1,q}(1))$, where 
$$c_r (i_1, \ldots, i_n)=\prod_{\stackrel{1 \leq k <r} {k \neq
\sum_{j=1}^s i_j }} k(r-k).$$
\begin{remark} \label{remarksing}
{\em Note that every singular vector 
(\ref{singular}) has form $L(-1)^{m+1}+\ldots$, where dots represent
lower degree terms (with respect to the universal enveloping algebra
grading).}
\end{remark}

\section{Vertex operator algebra $L(1,0)$}
\subsection{Zhu's algebra and intertwining operators}
We will use the definition of vertex operator algebra and
modules as stated in \cite{FHL} or \cite{FLM}.
Let $L(1,0)=M(1,0)/ \langle L(-1){\bf 1} \rangle$ be a 
simple vertex operator algebra associated to 
irreducible representation of the Virasoro algebra
(cf. \cite{FZ}, \cite{W}).

It is known that to every vertex
operator algebra $V$, we can associate Zhu's associative algebra $A(V)$ 
(cf. \cite{FZ} and \cite{Z}). 
In the special case $V=L(1,0)$, we know (see \cite{FZ}, \cite{W}) 
that $A(V) \cong {\mathbb C}[y]$,
where $y=[L(-2)-L(-1)].$
We have chosen the multiplication in 
$A(V)$ as in \cite{W} (which is slightly different
then the one in \cite{FZ}), 
$$a * b ={\rm Res}_{x} Y(a,x) \frac{(1-x)^{{\rm deg}(a)}}{x}b,$$
where $a, b \in A(V)$. 

By using standard techniques (see \cite{FZ}, \cite{W})
we have the following.
\begin{proposition}
Every irreducible module for the
vertex operator algebra $L(1,0)$ is isomorphic 
to $L(1,h)$ , for some $h \in {\mathbb C}$.
\end{proposition}
{\em Proof:} According to \cite{Z}, 
there is a one--to--one equivalence between 
equivalence classes of ${\mathbb N}$--gradable 
irreducible $L(1,0)$--modules and irreducible ${\mathbb C}[y]$--modules.
Every irreducible $L(1,0)$--module is a ${\rm Vir}$--module. Any such
module is $\mathbb{N}$--gradable and isomorphic to $L(1,h)$ for some
$h \in \mathbb{C}$. 
On the other hand
every finite dimensional irreducible ${\mathbb C}[y]$--module 
is one dimensional so the proof follows.
\epf

Since the notion of intertwining operator is more
subtle we include here the original definition 
\cite{FHL}.
\begin{definition} \label{intertwining}
{\em Let $W_1,W_2$ and $W_3$ be a triple of modules for vertex
operator algebra $V$. A mapping
$${\cal Y} \mapsto W_1 \otimes W_2 \rightarrow W_3\{x \},$$
is called an intertwining operator of type ${W_3 \choose W_1 \ W_2}$,
if it satisfies the following properties

\begin{enumerate}

\item The {\it truncation} property: For any $w_i \in W_i$, $i=1,2$,  
$$(w_1)_n w_2=0,$$
for $n$ large enough.

\item The {\it $L(-1)$-derivative property}: For any $v\in V$, 
$$\mathcal{Y}(L(-1) w_1, x)=\frac{d}{d x}
\mathcal{Y}(w_1, x),$$

\item The {\it Jacobi identity}: In ${ \rm Hom}(W_{1}\otimes W_{2}, W_{3})
\{x_{0},
x_{1},
x_{2}\}$, we have 
\bea
\lefteqn{x_0^{-1} \delta \left ( \frac {x_1-x_2}{x_0} \right ) 
Y(u, x_1) \mathcal{Y}(w_1, x_2)} \nn
&& -x_0^{-1} \delta \left ( \frac
{x_2-x_1}{-x_0} \right ) \mathcal{Y}(w_1, x_2) 
Y(u, x_1)
\nn
&&= x_2^{-1} \delta \left ( \frac {x_1-x_0}{x_2} \right
)\mathcal{Y}(Y(u, x_0)w_1, x_2)
\eea
for $u\in V$ and $w_1 \in W_1$.
\end{enumerate}
}
\end{definition}

We denote the space of all intertwining operators of the type 
${W_3 \choose W_1 \  W_2}$
by $I \ { W_3 \choose W_1 \  W_2}$.
The dimension of the space of intertwining operators (also
known as ``{\it fusion rule''} )of the type
${ W_3 \choose W_1 \ W_2}$ we denote
by ${\cal N}_{W_1,W_2}^{W_3}$.

Our goal is to 
find the fusion rules for the degenerate minimal models, i.e., 
$${\rm dim} \ I { L(1,\frac{r^2}{4}) \choose L(1,\frac{p^2}{4}) \ L(1,\frac{q^2}{4}) }.$$
Since our modules are irreducible we want to introduce Frenkel-Zhu's 
formula which gives us (roughly) a prescription for calculating 
fusion rules.
It is not hard to see, by using the Jacobi identity, that the space
$I{ L(1,\frac{r^2}{4}) \choose L(1,\frac{p^2}{4}) \ L(1,\frac{q^2}{4}) }$
is at most one dimensional.

Now for every module $M$, we associate an $A(V)$--bimodule
$A(M):=M/O(M)$ (cf. \cite{FZ}), where
$O(M)$ is spanned by the elements of the
form
$${\rm Res}_{x} Y(u,x) \frac{(1-x)^{deg(a)}}{x^2}v,$$
$u \in V$, $v \in M$.
In the case $M=M(c,h)$, 
\begin{equation} \label{bimodule}
O(M(c,h))= \{(L(-n-3)-2L(-n-2)+L(-1))v, n \geq 0, v \in M(c,h) \}.
\end{equation}
If we let
$$y=[L(-2)-L(-1)], \ \ x=[L(-2)-2L(-1)+L(0)],$$
then from the formulas
$$[L(-n)v]=[(ny-x+{\rm wt}(v))v],$$
and $$[x,y]w=0 \ {\rm mod} \ O(M(c,h)),$$
($[x,y]=xy-yx$) it follows that
$$A(M(c,h)) \cong {\mathbb C}[x,y],$$
as a $\mathbb{C}[y]$--bimodule (cf. \cite{L2}) ,
where the lowest weight vector is identified
with $1 \in \mathbb{C}[x,y]$ and the
actions of are
$$y*p(x,y)=xp(x,y), \ p(x,y)*y=yp(x,y),$$
for every $p(x,y) \in \mathbb{C}[x,y]$.

The Frenkel-Zhu's formula (\cite{FZ})
states that it is possible to calculate the dimension of the space
${M_{3}\choose M_{1} \ M_{2}}$ 
by knowing $A(V)$, $A(M_1)$, $M_2(0)$ and $M_3(0)$.
Instead of giving the original statement 
from \cite{FZ}, we state the following refinement 
obtained in \cite{L1}-\cite{L2}:
\begin{theorem} \label{lll}
Let $M_1$, $M_2$ and $M_3$ be lowest weight $V$--modules.
Suppose that $M_{2}$ and 
$M'_3$ are generalized Verma $V$--modules (see Section 3.2). Then
we have
$$\mathcal{N}_{M_{1} M_{2}}^{M_{3}} = \rm{dim} \ \rm{Hom}_{A(V)} 
(A(M_{1})\otimes _{A(V)}M_{2}(0), M_{3}(0)),$$
where $M_{i}(0)$, $i=1, 2, 3$,
is the ``top'' level of $M_{i}$, respectively,
equipped with the $A(V)$-module structures.
\end{theorem}

This theorem is not so useful as it stands.
On the other hand its proof is important. Hence it will be necessary
to understand a little bit deeper assumptions on $M_2$ and $M_3$ in our
situation.
For warm up let us
start with the ``easy--half'' of the Frenkel-Zhu's formula which says:
\begin{lemma} \label{upper}
Let $M_3$ be an irreducible lowest weight $V$--module. Then
$$ \mathcal{N}_{M_{1} M_{2}}^{M_{3}} \leq \rm{dim} \ \rm{Hom}_{A(V)} 
(A(M_{1})\otimes _{A(V)}M_{2}(0), M_{3}(0)).$$
\end{lemma}

Define an infinite dimensional Lie algebra ${\cal L}$ spanned
by $$L(-n-2)-2L(-n-1)+L(-n),$$
for $n \geq 1$.
In the case of minimal models--which is the most interesting 
case--the homology groups $H_q({\mathcal L},L(c,h))$
where calculated in \cite{FF2}.
For the Verma modules 
the $0$-th homology, $H_0({\cal L}, M(1,h))$ with the coefficients
in the Verma modules is isomorphic to ${\mathbb C}[x,y]$ as an
$A(L(1,0))$--bimodule (cf. \cite{W}).

The following result is an application
of a more general theory \cite{FF1}. 
\begin{theorem}
We have
\begin{enumerate}
\item[(a)]
$${\rm H}_0\left({\cal L}, L\left(1,\frac{m^2}{4}\right)\right),$$
is infinite--dimensional.
\item[(b)]
${\rm H}_0({\cal L}, L(1,\frac{m^2}{4}))$ is finitely generated as
a (left) $A(L(1,0))$--module.
\item[(c)]
$${\rm Ext}^1_{Vir, {\mathcal O}} 
\left(L(1,\frac{m^2}{4}),L(1,\frac{n^2}{4})\right) =\left\{\begin{array}{c}
{\mathbb C} \ \ {\rm if} \ \ |m-n|=2 \\
0  \ \ {\rm otherwise} \end{array} \right. $$
where ${\rm Ext}^1_{Vir, {\mathcal O}}$ stands for the relative ${\rm Ext}$ with respect
to the one--dimensional abelian subalgebra generated by $L(0)$.
\end{enumerate}
\end{theorem}
{\em Proof:}
a) Since the maximal submodule of $M(1,\frac{m^2}{4})$ is 
generated by one vector, in the projection (or homology) 
$A(L(1,\frac{m^2}{4}))$ is isomorphic to
$\frac{{\mathbb C}[x,y]}{I}$, where $I$ is a cyclic submodule
(with respect to the left and right actions) 
generated by some polynomial $p(x,y)$ which is a projection
of $v_{1,m}$ in ${\mathbb C}[x,y]$ .
It is clear that this space is infinite dimensional. \\
b)  Note first that $[L(-1)v]=(y-x+deg(v))[v]$. By using
Remark \ref{remarksing} it follows that
$$[v_{sing}]=p(x,y)=\prod_{i=1}^{m+1}(x-y+i)+q(x,y).$$
where $\mbox{deg}(q) < (m+1)$. Thus, the pure monomials 
in $p(x,y)$ with the highest powers 
are $x^{m+1}$  and $y^{m+1}$.
Since, $I$ is spanned by $p(x,y){\mathbb C}[x]$, here
we consider only the left action,   it follows that
$\frac{{\mathbb C}[x,y]}{I}$ is finitely generated. The similar
argument holds for the right action. \\
c) The idea is the same as in \cite{FF1}. The result is however different.
It is known that
$${\rm Ext}^*_{Vir, {\mathcal O}}(M,N) \cong H^*(Vir,{\mathcal O},Hom(M,N)).$$
Therefore
$$ H^*(Vir,{\mathcal O},Hom(M,N)) \cong Tor_*^{Vir,\mathcal{O}}(N^*,M),$$
where $N^*$ is the dual module.
Hence we can compute our cohomology by using the tensor product
of complexes 
$$M(1,\frac{(m+2)^2}{4}) \longrightarrow  M(1,\frac{m^2}{4}),$$
$$M(1,\frac{(n+2)^2}{4})^{opp} \longrightarrow  M(1,\frac{n^2}{4})^{opp},$$
where $M(c,h)^{opp}$ is the opposite Verma module (cf. \cite{FF1}-\cite{FF2}). 
The corresponding spectral sequence $E_2^{p,q}$ collapses at the second term
Therefore 
$$Tor_1^{Vir,\mathcal{O}}(L(1,\frac{n^2}{4})^*,L(1, \frac{m^2}{4}) \cong E_2^{1,0}
\cong \mathbb{C} \ {\rm or} \ 0,$$
where non--trivial homology occurs 
only if the Verma module $M(1,\frac{m^2}{4})$ embeds inside $M(1,\frac{n^2}{4})$ as the maximal submodule or vice--versa.
This happens if and only if $|n-m|=2$. Therefore we have the proof \footnote{It is crucial to
notice that our cohomology is relative one, otherwise our extension are not controllable
inside category $\mathcal{O}$. Such (non--relative) extensions are studied in \cite{M}}.
The corresponding short--exact sequences are clearly,
\begin{equation} \label{ext}
0 \rightarrow L(1,\frac{(m+2)^2}{4}) \rightarrow M(1, \frac{m^2}{4})/M(1,\frac{(m+4)^2}{4})
\rightarrow L(1,\frac{m^2}{4}) \rightarrow 0,
\end{equation}
and the one obtained from \ref{ext} by applying (exact) functor $( \ )'$ taking 
modules to the corresponding contragradient modules.
\epf

For every $m,n \in {\mathbb N}$ (we exclude the case $mn=0$), 
fix a multiset $J_{m,n}=\{m+n,m+n-2, \ldots,m-n \}$.
Let ${\cal F}_{\lambda, \mu}$ be a ``density'' module 
for the Virasoro algebra. ${\cal F}_{\lambda,\mu}$
is spanned by $w_r$, $r \in {\mathbb Z}$  and the action is given by
$$L_n. w_r=(\mu+r+\lambda(m+1))w_{r-n}.$$
In \cite{FF1} the projection formula
for the singular vectors (considered as an
element of the enveloping algebra) on ${\cal F}_{\lambda,\mu}$
(more precisely $w_0$) was found.
We want to relate the projection of the singular
vectors on ${\cal F}_{\lambda, \mu}$ with the projection
inside $A(M(1,\frac{m^2}{4}))\otimes_{C[y]} L(1, \frac{n^2}{4})$).
It is easy to see that
\bea \label{projectionam}
&& [L(-j_1)\ldots L(-j_k)v_{m^2/4}]= \nn
&& \prod_{r=1}^k (j_r \frac{n^2}{4}-y+\beta(r,k)).[v_{m^2/4}]= \nn
&& \prod_{r=1}^k(j_r \frac{n^2}{4}-x+\beta(r,k))v_{m^2/4}
\eea
where $v_{m^2/4}$ is the lowest weight vector and
$$\beta(r,k)=j_{r+1}+\ldots+j_k+\frac{m^2}{4}.$$
But the last factor in
(\ref{projectionam}) is the same as
the $P(j_1,..,j_k)$  
where
$$L(-j_1)\ldots L(-j_k).w_0=P(j_1,...,j_k)w_{j_1+...+j_k},$$
and the projection is in ${\cal F}_{\lambda,\mu}$ 
for $\lambda=-\frac{n^2}{4}$ and $\mu=\frac{n^2}{4}+\frac{m^2}{4}-x$.

In the remarkable paper \cite{FF2}, 
projection formulas for all singular vectors on the density modules
were found. In the slightly different notation, for the singular vectors we consider, 
these formulas appeared in \cite{K}.
The result is
\begin{equation} \label{ffkent}
v_{1,m+1}.w_0=\prod_{i \in J_{m,n}}(x-\frac{i^2}{4})w_{m+1},
\end{equation}
up to a multiplicative constant.
  
Now, by using (\ref{ffkent}) fact and 
the discussion above (cf. \cite{W})
we obtain
\begin{lemma} \label{Kent}
As a $A(L(1,0)$--module 
$A(L(1,\frac{m^2}{4})) \otimes_{A(L(1,0))} L(1,\frac{n^2}{4})(0)$
is isomorphic to
$\frac{{\mathbb C}[x]}{<\prod_{i \in J_{m,n}}(x-i^2/4)>}.$
\end{lemma}

If $n \leq m$ notice 
that as an $A(L(1,0))$-- module
\begin{equation} \label{avdecoposition}
A(L(1,\frac{m^2}{4})) \otimes_{A(L(1,0))} L(1,\frac{n^2}{4})(0)
\cong \bigoplus_{i \in J_{m,n}} {\mathbb C}v_i ,
\end{equation}
where $v_i$ is an irreducible $A(L(1,0))$--module
such that $y.v_i=i^2/4 v_i$.
But if $m < n$, then we have two-dimensional
submodule in the 
above decomposition (and this module is {\em not} completely
reducible). 
Thus,  (\ref{avdecoposition}) is {\em not} symmetric if we switch
$m$ and $n$.

The similar failure was already noticed  in \cite{L1}.
Anyhow, by using Lemma \ref{Kent} and Lemma \ref{upper}
we obtain 
\begin{proposition} \label{upperprop}
Let $L(1, \frac{m^2}{4})$, $L(1, \frac{n^2}{4})$ and $M$ an irreducible
$L(1,0)$--modules.
Then we have the following upper bounds 
\begin{eqnarray} \label{upper1}
&& {\rm dim} \ I {M \choose L(1,\frac{m^2}{4}) \  L(1,\frac{n^2}{4}) } \leq
\left\{ \begin{array}{c} 1 \  {\rm if} \ M \cong
L\left(1,\frac{r^2}{4}\right) \  {\rm for} \ r \in J_{m,n}, 
\\ 0 \ {\rm otherwise} \end{array} \right.
\end{eqnarray}
where $J_{m,n}=\{m+n, \ldots, m-n \}$.
\end{proposition}
\epf

Now, we shall show that the  equality
holds in the equation (\ref{upper1}). We will provide 
two different proofs. One which uses the properties of
Verma modules and the other which uses free field realization
of the modules $L(1,\frac{m^2}{4})$.

\subsection{Lie algebra $g(V)$}
Let $V$ be a vetex operator algebra. Let $\hat{V}=V \otimes \mathbb{C}[t,t^{-1}]$,
$d=L(-1) \otimes 1 + 1 \otimes \frac{d}{dt}$ and $g(V)=V / dV$.
It has been noticed by several authors that the space $g(V)$ 
has a Lie algebra structure if we let
$$[a(m),b(n)]=\sum_{i=0}^{\infty}{m \choose i} (a_i b)(m+n-i).$$
If we define the grading 
with ${\rm deg}a(m)=n-m-1$, where $a \in V_{(n)}$, then
we have the corresponding triangular decomposition
$g(V)=g(V)_- \oplus g(V)_0 \oplus g(V)_+$.
Let $U$ be any $g(V)_0$--module.
We let (as in \cite{L2})
$$F(U)={\rm Ind}_{U(g(V)_+ \oplus g(V)_0)}^{U(g)} U,$$
such that $g(V)_+$ acts as zero.
We define also the quotient $\bar{F}(U)=F(U)/J(U)$ (the so--called generalized Verma module \cite{L2}), 
where $J(U)$ is the
intersection of all kernels of all $g(V)$--homomorphisms from $F(U)$ to
weak modules. Now, the assumption in Theorem \ref{lll} on $M_2$
and $M_3'$ means that
$M_2 \cong \bar{F}(M_2(0))$ and $M_3=\bar{F}(M_3^*(0))'$.

In \cite{L1}--\cite{L2} it was shown that
every $A(V)$ homomorphism 
from $A(W_1) \otimes_{A(V)} W_2(0)$ to $W_3(0)$
does not necessary lead to an
intertwining operator of the type 
${W_3 \choose W_1 \ W_2}$ but rather
to ${F(W_3(0)^*)' \choose W_1 \ F(W_2(0)) }$ (actually $F(W_2(0))$
might be replaced by $\bar{F}(W_2(0))$).

In the case when $V$ is rational
$\bar{F}(W_2(0)) \cong W_2$ and $\bar{F}(W_3(0)^*)' \cong W_3$
(\cite {L2}). 
But if the vertex operator algebra $V$ 
is not rational, the main difficulty is that the generalized 
Verma module $\bar{F}(W_2(0))$ may not be isomorphic
to  $W_2$ (let alone $F(W_2(0))$ !) (cf. \cite{L2}). 
Also, the spaces $F(U)$ and $\bar{F}(U)$ are extremely
difficult to analyze explicitly.
Still, because we are dealing with a particular example, Virasoro vertex operator algebra, we
can make use of singular vectors and Verma modules to
simplify the whole construction.

Let $V=L(1,0)$.
Pick $\omega =L(-2){\bf 1} \in L(1,0)$.  
Then, inside $g(L(1,0))$, we have 
$$[\omega(m+1),\omega(n+1)]=(m-n)\omega(m+n+1)+
\delta_{m+n,0} \frac{m^3-m}{12},$$
i.e., these operators close the Virasoro algebra.
From the construction of $F(U)$ it is clear that  
$U(Vir_-) \otimes U \hookrightarrow U(g(V)_-) \otimes U \cong F(U)$. 
In particular
$M(1,h) \hookrightarrow F(M(1,h)(0))$.


\subsection{The fusion rules computations}

Assume first that
\begin{equation} \label{mn}
m \leq n.
\end{equation}
First we replace the ``big'' space $F(M(1,h))$  
with the smaller Verma module for the Virasoro algebra (we have
seen already that the latter is a subspace inside $F(M(1,h))$).

Now, let us pick a non--trivial $A(L(1,0))$ homomorphism from
$A(L(1,\frac{m^2}{4})) \otimes_{A(L(1,0))} L(1,\frac{n^2}{4})(0)$ to
$L(1,\frac{r^2}{4})(0)$. Also let
$T=L(1,\frac{m^2}{4}) \otimes 
\mathbb{C}[t,t^{-1}] \otimes M(1,\frac{n^2}{4}$ be a $g(L(1,0))$--module 
as in \cite{L2}. Then the construction in 
\cite{L2} implies that there is a bilinear pairing between
$T$ and $M(1,\frac{r^2}{4}) \hookrightarrow F(M(1,\frac{r^2}{4})(0)^*)$.
This implies (again by applying Li's construction in the proof of Theorem 2.11 in \cite{L2}) that the corresponding
intertwining operator {\em lands} in
$M(1,\frac{r^2}{4})'$, i.e., it is of the type
${M(1,\frac{r^2}{4})' \choose L(1,\frac{m^2}{4}) \ M(1,\frac{n^2}{4})}$.
Here $M(1, \frac{r^2}{4})'$ is the contragradient Verma module
(cf. \cite{FF2}). The contragradient module
$M(1, \frac{r^2}{4})'$ is {\em not} of the lowest weight type 
(because $M(1, \frac{r^2}{4})$ is reducible). 
In particular, if $v '$ is the lowest weight vector 
$$U(Vir)v ' \cong L(1, \frac{r^2}{4}),$$ 
i.e. we can ``paste'' the whole irreducible module by acting on
the lowest weight subspace, but not the whole module $M(1,\frac{m^2}{4})'$. 
Now, the question is \\
{\em How to descend from $M(1,\frac{m^2}{4})'$ to 
$L(1,\frac{m^2}{4})$ ?} \\
Here is the proof.
We have either $n \leq r$ or $r < n$.
For each of these two cases we consider 
\begin{equation} \label{prvislucaj}
I{M(1,\frac{r^2}{4})' \choose L(1,\frac{m^2}{4}) \ M(1,\frac{n^2}{4})},
\end{equation}
or
\begin{equation} \label{drugislucaj}
I{M(1,\frac{n^2}{4})' \choose L(1,\frac{m^2}{4}) \ M(1,\frac{r^2}{4})},
\end{equation}

respectively.
Notice that these two spaces are isomorphic because
of $I { M_3 \choose M_1 \ M_2} \cong I{ M'_2 \choose M_1 \ M'_3}$.
Suppose that $n \leq r$. 

Now the aim is to construct 
intertwining operator 
of the type ${ M(1,\frac{r^2}{4})' \choose L(1,\frac{m^2}{4})
\ L(1,\frac{n^2}{4})}.$
Therefore if we can check
\begin{equation} \label{sss}
\langle w'_3,{\mathcal Y}(w_1,x)w \rangle =0,
\end{equation}
for every $w \in M(1,\frac{(m+2)^2}{4}) \hookrightarrow
M(1,\frac{m^2}{4})$, $w'_3 \in M(1,\frac{r^2}{4})''=M(1,\frac{r^2}{4})$ and
$w_1 \in L(1,\frac{n^2}{4})$, then by defining
$\bar{\mathcal Y}(w_1,x)[w_2]:=\bar{\mathcal Y}(w_1,x)w_2$
where $[w_2] \in M(1,\frac{m^2}{4})/M(1,\frac{(m+2)^2}{4})$,
we obtain a (well--defined) non--trivial intertwining operator
of the type ${ M(1,\frac{r^2}{4})' \choose L(1,\frac{m^2}{4})
\ L(1,\frac{n^2}{4})}$.

Let us check that (\ref{sss}) holds. First of all, because
of the Jacobi identity and the fact that $M(1,\frac{r^2}{4})$ 
is lowest weight module, it is enough to show that
\begin{equation} \label{uh}
\langle w'_3,{\mathcal Y}(w_1,x) v_{sing} \rangle =0
\end{equation}
where $w'_3 \in M(1,\frac{r^2}{4})''(0)=M(1,\frac{r^2}{4})(0)$
is the lowest weight vector and $v_{sing}$ is the singular vector 
that generates the maximal submodule of $M(1,\frac{m^2}{4})$.
\bea \label{oko}
&& \langle w'_3, {\mathcal Y}(w_1,x)L(-j_1) \ldots L(-j_k) w \rangle= \nn
&& \prod_{i=1}^k -(x^{-{j_i}+1}\partial_x + (1-j_i)x^{-j_i}\frac{n^2}{4}) \langle w'_3, 
{\mathcal Y}(w_1,x) w \rangle= \nn
&& \prod_{i=1}^k -(x^{-{j_i}+1}\partial_x + (1-j_i)x^{-j_i}\frac{m^2}{4})
C x^{\frac{r^2}{4}-\frac{m^2}{4}-\frac{n^2}{4}}=\nn
&& (-1)^{\sum_{i} j_i}\prod_{i=1}^k
(\frac{r^2}{4}-\frac{m^2}{4}-\frac{n^2}{4}-
\sum_{s=i+1}^{k} j_s+(1-j_i)\frac{m^2}{4})C 
x^{\frac{r^2}{4}-\frac{m^2}{4}-\frac{n^2}{4}}=\nn
&& C \prod_{i=1}^k (j_i \frac{m^2}{4}-\frac{r^2}{4}
+\sum_{s=i+1}^{k} j_s+\frac{n^2}{4})C
x^{\frac{r^2}{4}-\frac{m^2}{4}-\frac{n^2}{4}-\sum_{i} j_i},
\eea
where $C$ is a constant that depends on $\mathcal Y$ (we may assume
that $C$ is equal to $1$).
If we compare (\ref{oko}) with (\ref{projectionam}) we see  
that products appearing in both expressions are the same if we interchange $x$ 
with $\frac{r^2}{4}$ and $\frac{m^2}{4}$ with $\frac{n^2}{4}$. 
In other words the expression
$\langle w'_3,{\mathcal Y}(w_1,x) v_{sing} \rangle =0$
if and only if the corresponding projection inside
$A(L(1,\frac{n^2}{4}) \otimes_{A(L(1,0)}
A(L(1,\frac{m^2}{4}))$
is zero (notice that now $L(1,\frac{n^2}{4})$ and
$L(1,\frac{m^2}{4})$ changed positions).
We know that
$$A(L(1,\frac{n^2}{4}) \otimes_{A(L(1,0)} A(L(1,\frac{m^2}{4})) \cong
\frac{\mathbb{C}[x]}{\prod_{i \in J_{n,m}} \langle x-\frac{i^2}{4}
\rangle }.$$
Because of (\ref{mn}), $J_{n,m} \subset J_{n,m}$ (as multisets). Therefore  
$$\langle w'_3,{\mathcal Y}(w_1,x) v_{sing} \rangle =0$$
holds.
Thus we obtain 
a non--trivial  intertwining operator
$\bar{\cal Y}$ of the type
${ M(1,\frac{r^2}{4})' \choose L(1,\frac{m^2}{4}) \ L(1,\frac{n^2}{4})}.$
Now,
$$I { M(1,\frac{r^2}{4})' \choose L(1,\frac{m^2}{4}) \ L(1,\frac{n^2}{4})}
\cong I{ L(1,\frac{n^2}{4}) \choose L(1,\frac{m^2}{4}) \ M(1,\frac{r^2}{4})}.$$
Because of our initial assumption $n \leq r$, and $m-n \leq r \leq m+n$
it follows that $m-r \leq n \leq m+r$, therefore we can repeat the whole
procedure for  $M(1,\frac{r^2}{4})$ so we end up with a non--trivial
intertwining operator of the type
$${ L(1,\frac{n^2}{4}) \choose L(1,\frac{m^2}{4}) \ L(1,\frac{r^2}{4})}.$$
If $r<q$ then we pick the intertwining operator (\ref{drugislucaj})
and the same reasoning leads to a non--trivial intertwining operator
of the type 
$${ L(1,\frac{r^2}{4}) \choose L(1,\frac{m^2}{4}) \ L(1,\frac{n^2}{4})}.$$
This also follows from the duality property for the intertwining operators.
If we summarized everything we obtain
\begin{theorem} 
$${\rm dim} \ I {L(1,\frac{r^2}{4}) \choose L(1,\frac{m^2}{4}) \ 
L(1,\frac{n^2}{4})}=1$$
if and only if $r \in \{m+n,...,|m-n|\}$.
\end{theorem}

\begin{theorem}
Let ${\cal A}$ be a free Abelian group on the
set $\{ a(m) : m \in {\mathbb N} \}$
and
$$\times : {\cal A} \times {\cal A} \rightarrow {\cal A}$$ 
a binary operation defined 
by the formula
$$a(m) \times a(n)=\sum_{r \in {\mathbb N}\cup 0} {\cal N}_{L(1,\frac{m^2}{4})
L(1,\frac{n^2}{4})}^{L(1,\frac{r^2}{4})} a(r).$$
Then ${\cal A}$ is a commutative associative ring
with the multiplication
$$a(m) \times a(n)=a(m+n) + a(m+n-2)+ \ldots + a(|m-n|),$$ 
i.e. ${\cal A}$ is isomorphic to the
representation ring ${\cal R}ep({\mathfrak sl}(2,{\mathbb C}))$.
\end{theorem}
\epf

\begin{remark}
{\em In general if $M$ is any $L(1,0)$--module 
and
$${\cal Y} \in I{M \choose L(1,\frac{m^2}{4}) \ L(1,\frac{n^2}{4})},$$
then $M$ is not necessary completely reducible.
Also, note that we excluded the case $mn=0$. If $m$ or $n$ are equal to
zero then we deal with intertwining operators among two irreducible modules and
vertex operator algebras, which are well known.

Another interesting fact is that in the case (\ref{mn}) the 
module $A(L(1,\frac{m^2}{4})) \otimes_{A(L(1,0))}
L(1,\frac{n^2}{4})(0)$
is not completely reducible. This fact was exploited in \cite{M}
where we study {\em logarithmic} intertwining operators.}
\end{remark} 

Note that in our proof we actually analyzed more
carefully the failure of Frenkel-Zhu's formula.
One should not expect to apply our procedure 
in the more general setting, because our Virasoro
vertex operator algebra has a quite simple structure.
Certainly it would be interesting to
study a class of vertex operator algebra for which
\begin{equation} \label{symmetry}
A(W_1) \otimes_{A(V)} W_2(0) \cong A(W_2) \otimes_{A(V)} W_1(0),
\end{equation}
for any choice of irreducible modules $W_1$ and $W_2$.
Then we hope that for this 
class of vertex algebras some version of Frenkel-Zhu's formula
indeed apply. Assumption  (\ref{symmetry}) turns out to be 
very natural since
\begin{equation} \label{ff}
I {W_3 \choose W_1 \ W_2} \cong I {W_3 \choose W_2 \ W_1}.
\end{equation}

\section{Construction of all intertwining operators
for the family ${\cal F}_1$}

\subsection{$V_L$ vertex operator algebra  and its irreducible modules}

Let $L$ be a rank one even lattice with a generator $\beta$ normalized
such that $\langle \beta,\beta \rangle =1$  and let
$\alpha=\sqrt{2} \beta$. Thus $\langle \alpha,\alpha
\rangle=2$.
As in \cite{FLM}, \cite{DL} we define $V_L$ as a vector space
$$V_L =M(1) \otimes {\mathbb C}[L],$$
where $M(1)$ is the level one irreducible
module for  Heisenberg algebra $\hat{h}_{Z}$ associated
to one--dimensional abelian algebra $h=L \otimes_{\mathbb Z} {\mathbb C}$ and
${\mathbb C}[L]$ is the group algebra of $L$ with a 
generator $e^{\alpha}$. Put
$\omega=\frac{1}{2}\beta(-1)^2$. Then $V_L$ is a vertex
operator algebra (see \cite{FLM}) with the Virasoro element $\omega$.
We have  a decomposition
$$V_L = \bigoplus_{m \in {\mathbb Z}} M(1) \otimes e^{m \alpha}.$$
Let $L^o$ be a dual lattice,  $L^o/L \cong {\mathbb Z}/2{\mathbb Z}$.
Then (as in \cite{DL}), for a nontrivial coset
representative, we obtain an irreducible $V_L$--module
$V_{L+1/2}$, which
can be decomposed as
$$V_{L+1/2} = \bigoplus_{m \in {\mathbb Z}} M(1) \oplus e^{m \alpha + 1/2
\alpha}.$$
Moreover, $V_{L+1/2}$, $V_L$  is (up to equivalence) complete
list of irreducible $V_L$--modules.
Furthermore, one can equip the
space $W=V_L \oplus V_{L+1/2}$ (as in \cite{DL}) with
the structure of the generalized vertex operator algebra.
We will neglect this fact in our considerations.

For every module $W$ for the Virasoro algebra
on which $L(0)$ acts semisimple we define a formal character
(or a $q$-graded dimension) by 
$$ch_{q}(W)=\sum_{ n \in {\rm Spec}L(0)} {\rm dim}(W_{n}) q^n.$$
 
From the Proposition \ref{sequence} 
it follows that
$$ch_{q}(L(1,\frac{m^2}{4}))=\frac{q^{\frac{m^2}{4}}-q^{\frac{(m+2)^2}{4}}}{q^{-1/24}\eta(q)}.$$
Then it is not hard to obtain
\begin{eqnarray} \label{characters}
&& ch_q(V_L)=\sum_{n \geq 0} (2n+1) ch_q(L(1,n^2)) \nn
&& ch_q(V_{L+1/2})=\sum_{n \geq 0} (2n+2) ch_q
\left(L\left(1,\frac{(2n+1)^2}{4} \right) \right).
\end{eqnarray}

Consider the vectors
$$x=e^{\alpha}, \ y=e^{-\alpha}, \ h=\alpha(-1) \iota(0),$$
which span $(V_L)_1$. These vectors span a Lie
algebra isomorphic to ${\mathfrak sl}(2,{\mathbb C})$.
$x_0$, $y_0$ and $h_0$  as act derivatives on $W$.
The following result was
obtained in \cite{DG}.

\begin{proposition} \label{evendec}
As $(L(1,0), \mathfrak{sl}_2)$--module
$$V_L \cong \bigoplus_{m \geq 0} L(1,m^2) \otimes V(2m),$$
where $V(2m)$ is an irreducible $2m+1$ dimensional $\mathfrak{sl}_2$--module.
\end{proposition}
\epf

The proof uses the result from \cite{DLM}, \cite{DM1} about the decomposition
of the vertex operator algebra $V$ with respect to
a ``dual'' pair $(V^G, G)$ where $G = Aut(G) $ is a compact (or finite) group
and $V^G$ is a $G$--stable subvertex operator algebra.
This can be modified when instead of group $G$ 
we work with the Lie algebra.

Since $V_{L+1/2}$ is a module for the pair $(V^{\mathfrak{sl}_2}_L, 
\mathfrak{sl}_2)$ then
by using (\ref{characters}) we derive
\begin{equation} \label{odddec}
V_{L+1/2} \cong \bigoplus_{m \geq 0}  L(1, \frac{(2m+1)^2}{4}) \otimes
V(2m+1),
\end{equation}
where $V(2m+1)$ is a $2m+2$--dimensional  $\mathfrak{sl}_2$--module.
It easy to see that 
$V(2m+1)$ is irreducible $\mathfrak{sl}_2$--module.

\begin{remark} {\em Note that $V^{sl_2}$ ($sl_2$--stable vertex operator
algebra) is exactly $V^{G}$ where $G \cong SO(3)$ is a (full) 
group of automorphisms of $V_L$. 
It is well known that every irreducible representation can be
obtain as a representation of $SL(2,{\mathbb C})$, since
$PSL(2,{\mathbb C}) \cong SO(3)$. 
In particular
every such finite-dimesional representation is odd--dimensional.}
\end{remark}

Since, 
$V_{L+1/2}$ is an irreducible $V_L$--module
we have the  
Jacobi  identity 
\bea \label{jacobim}
&& x_0^{-1} \delta \left ( \frac {x_1-x_2}{x_0} \right ) 
Y(u, x_1){Y}(v, x_2)w \nn
&& -x_0^{-1} \delta \left ( \frac
{x_2-x_1}{-x_0} \right ){Y}(v, x_2) 
Y(u, x_1)w
\nn
&&= x_2^{-1} \delta \left ( \frac {x_1-x_0}{x_2} \right
){Y}(Y(u, x_0)v, x_2)w,
\eea
for every  $u \in V_{L}$, $v \in V_{L+1/2}$ and $w \in W$. 
Also, for 
$${\cal Y} \in I {V_L \choose V_{L+1/2} \ V_{L+1/2}},$$  
we have
\bea \label{jacobig}
&& x_0^{-1} \delta \left ( \frac {x_1-x_2}{x_0} \right ) 
Y(u, x_1) {\cal Y}(v, x_2)w \nn
&& -x_0^{-1} \delta \left ( \frac
{x_2-x_1}{-x_0} \right ){\cal Y}(v, x_2) 
Y(u, x_1)w
\nn
&&= x_2^{-1} \delta \left ( \frac {x_1-x_0}{x_2} \right
)\mathcal{Y}(Y(u, x_0)v, x_2)w.
\eea

\begin{remark} \label{remaint}
{\em Note that
$W$ can not be equipped with a vertex operator superalgebra structure.
If $u,v \in V_{L+1/2}$ then we do not get Jacobi identity in
the form (\ref{jacobim}) or (\ref{jacobig}), but rather generalized 
identity where the delta function is suitably multiplied with the
terms of the type $\left(\frac {x_1-x_0}{x_2} \right)^{1/2}$.
Studying this (generalized) Jacobi identity is useful
for studying convergence and the extension properties
for the intertwining operators (cf.  \cite{H1}).}
\end{remark}

\subsection{Intertwining operators for the family ${\cal F}_1$.}

Let $V(i)$, $i \in {\mathbb N}$  be an irreducible $sl_2$--module considered as a 
subspace of $W$ which corresponds to the 
decompositions in Proposition \ref{evendec} and (\ref{odddec}).
Fix a positive integer $j$. 
We introduce a basis $u_j(m)$, $m \in \{j, j-2,\ldots,-j \}$
for $V(j)$, such that the following relations are satisfied,
\begin{eqnarray}
&& h.u_j(m)=mu_j(m) \nn
&& x.u_j(m)=\frac{\sqrt{(j+m+2)(j-m)}}{2}u_j(m+2) \nn
&& y.u_j(m)=\frac{\sqrt{(j+m)(j-m+2)}}{2} u_j(m-2),
\end{eqnarray}
where $u_j(k)=0$ for $k \notin \{j, \ldots, -j \}$.
Also, we choose a dual basis $u^*_j(m)$ for
$V(j)^*$ such that
$<u^*_j(m),u_j(n)>=\delta_{m,n}$. Define
$<g.u^*,v>=-<u^*,g.v>$. Then $V(j)^*$ became
a $sl_2$--module and an isomorphism from
$V(j)$ to $V(j)^*$ is given by $\mu(u_j(m))=(-1)^{j-m}u_j^*(-m).$
By using this identification, 
for $j_1, j_2, j_3 \in {\mathbb N}$ and $-j_i \leq m_i \leq j_i$,
$i=1,2,3$, 
we introduce real numbers (Clebsch--Gordan coefficients) 
$\left( \begin{array}{ccc} j_1 & j_2 & j_3 \\ m_1 & m_2 & m_3 
\end{array} \right),$
such that
\begin{equation} \label{clebsch}
u_{j_1}(m_1) \otimes u_{j_2}(m_2)=\sum_{j_3=|j_1-j_2|}^{j_3=j_1+j_2}
\left( \begin{array}{ccc} j_1 & j_2 & j_3 \\ m_1 & m_2 & m_1+m_2 
\end{array} \right) u_{j_3}(m_1+m_2).
\end{equation}

First we need an auxiliary result which is slightly modified result
from  \cite{DM1} and \cite{DG}.
\begin{proposition} \label{independence}
Suppose that $V$ is a vertex operator algebra and $W_1$, $W_2$ and
$W_3$ three irreducible $V$--modules.
Let $v_i \in W_1 , w_i \in W_2$, 
$i=1,...,k$ be homogeneous elements  such that $v_i \neq 0$ and
$w_i$ are linearly independent.
Then
$$\sum_{i=1}^k {\cal Y}(v_i,x)w_i \neq 0.$$
\end{proposition}
\epf

Now. let us go back to our vertex operator algebra $V_L$.
Let ${\cal Y}$ be any intertwining operator of the type
\begin{equation} \label{triple}
{V_L \choose V_{L+1/2} \ V_{L+1/2}}, {V_{L+1/2} \choose V_{L}
\ V_{L+1/2}}  \ {\rm or} \ {V_L \choose V_{L} \ V_{L} }.
\end{equation}

By using the Proposition \ref{independence} the map
$${\cal Y}( \ ,x) \ : V(j_1) \otimes V(j_2) \rightarrow W\{x \}$$
is injective.
and for every $m_1, m_2$ and $j_1$, $j_2$ there is a $p \in {\mathbb C}$
such that
$$u_{j_1}(m_1)_p u_{j_2}(m_2)=\sum_{j_3=|j_1-j_2|}^{j_3=j_1+j_2} 
 k(j_1, j_2, j_3, m_1,m_2, m_1+m_2)   u_j(m_1+m_2),$$
where
$ k(j_1, j_2, j_3, m_1,m_2, m_1+m_2) $ is a (non--zero) multiple
of 
$$ \left( \begin{array}{ccc} j_1 & j_2 & j_3 \\ m_1 & m_2 & m_1+m_2
\end{array} \right).$$
(in the special case ${\cal Y}=Y$ this fact  was noticed in
\cite{DG}).

Now it is clear that if
$\left( \begin{array}{ccc} j_1 & j_2 & j_3 \\ m_1 & m_2 & m_1+m_2
\end{array} \right) \neq 0$, then
the $L(1,0)$--module generated by 
${\cal Y}( u_{j_1}(m_1), x)u_{j_2}(m_2)$ contains a copy
of $L(1, \frac{j_3^2}{4})$. 
Since $L(1,0)$ is contained in
$V_L$ and $L(1, \frac{m^2}{4})$ is an $L(1,0)$--module
then we obtain the following Jacobi identity
\bea \label{jacobih}
&& x_0^{-1} \delta \left ( \frac {x_1-x_2}{x_0} \right ) 
Y(u, x_1) {\cal Y}(v, x_2)w -x_0^{-1} \delta \left ( \frac
{x_2-x_1}{-x_0} \right ){\cal Y}(v, x_2) 
Y(u, x_1)w
\nn
&&= x_2^{-1} \delta \left ( \frac {x_1-x_0}{x_2} \right
)\mathcal{Y}(Y(u, x_0)v, x_2)w,
\eea
for $u \in L(1,0)$, $v \in L(1,\frac{j_1^2}{4})$ and $w \in L(1,
\frac{j_2^2}{4})$ (here $v$ and $w$ lie in Vir--submodules
generated by $u_{j_1}(m_1)$ and $u_{j_2}(m_2)$, respectively).

Now we can push down ${\cal Y}$ to 
$L(1,\frac{j_3^2}{4})$, which is generated by the vector $u_{j_3}(m_1+m_2)$,
since for every $j_1,j_2$ and $|j_1-j_2| \leq j_3 \leq j_1+j_2$
we can choose a pair $m_1$, $m_2$ and a ${\cal Y}$ of the
appropriate type (\ref{triple}) 
such that 
$$\left( \begin{array}{ccc} j_1 & j_2 & j_3 \\ m_1 & m_2 & m_1+m_2
\end{array} \right) \neq 0.$$
We obtain an intertwining operator of the type
$${L(1, \frac{j_3^2}{4}) \choose L(1, \frac{j_1^2}{4}) \ L(1,\frac{j_2^2}{4})},$$
and this is the end of the construction.

\section{Lie superalgebra $\mathfrak{osp}(1|2)$ and ${\rm
Rep}(\mathfrak{osp}(1|2))$}
The Lie superalgebra $\mathfrak{osp}(1|2)$ is a graded extension of
the finite--dimensional Lie algebra $\mathfrak{sl}(2,{\bf C})$. 
It has three even generators $x,y$ and $h$, and two odd generators
$\varphi$ and $\chi$, that satisfy:

$$[h,x]=2x, \ \ [h,y]=-2y, \ \ [x,y]=h,$$
$$[x,\chi]=\chi, \ \  [x,\varphi]=-\varphi, \ \ [y,\chi]=-\chi, \ \
[y,\varphi]=\varphi,$$
$$[h,\varphi]=-\varphi, \ \ [h,\chi]=\chi,$$
$$\{\chi,\varphi\}=2h, \{\chi,\chi\}=2x, \ \   \{\varphi,\varphi\}=2y .$$

Generators $\{x,y,h\}$ span a Lie algebra isomorphic to $\mathfrak{
sl}(2,\mathbb{C})$, and this fact
makes the representation theory of $\mathfrak{osp}(1|2)$ quite simple.
All irreducible $\mathfrak{osp}(1|2)$--modules can be constructed in the
following way. Fix a positive half integer $j$ ($2j \in {\bf N}$) and
a $4j+1$--dimensional vector
space $V(j)$ spanned by the vectors $\{v_{j},v_{j-1/2},...,v_{-j}\}$, with the
following actions: 
\begin{eqnarray}
&& x.v_{i}=\sqrt{[j-i][j+i+1]}v_{i+1}, \nn
&& y.v_{i}=\sqrt{[j+i][j-i+1]}v_{i-1}, \nn
&& h.v_{i}=2iv_{i}. 
\end{eqnarray}
If $2(i-j) \in {\bf Z}$ then we define
\begin{eqnarray}
&& \varphi.v_{i}=- \sqrt{j+i} v_{i-1/2}, \nn 
&& \chi.v_{i}=-\sqrt{j-i}v_{i+1/2},
\end{eqnarray}
otherwise
\begin{eqnarray}
&& \varphi.v_{i}=\sqrt{j-i+1/2} v_{i-1/2}, \nn 
&& \chi.v_{i}=-\sqrt{j+i+1/2}v_{i+1/2}.
\end{eqnarray}
In all these formulas $v_j=0$ if $j \notin \{j,j-\frac{1}{2},\ldots,-j \}$.
It is easy to see that each $V(j)$ is an irreducible $\mathfrak{osp}(1|2)$--module and
that every finite dimensional irreducible representation of $\mathfrak{osp}(1|2)$ is isomorphic to $V(j)$ for some $j \in {\mathbb N}/2$.

The representations with $j \in \mathbb{N}$ we call {\em even}, and the
representations with $j \in \mathbb{N}+\frac{1}{2}$ we call {\em odd}.
We extend this definition for an arbitrary element of
$V \in {\rm Rep}(\mathfrak{osp}(1|2))$. The corresponding decomposition is
$V=V_{{\rm even}}+V_{{\rm odd}}$.

It is a pleasant exercise to decompose the 
tensor product $V(i) \otimes V(j)$. The following result is
well--known:
\begin{equation} \label{osp12}
V(i) \otimes V(j) \cong  \bigoplus_{k=|i-j|, k \in {\bf N}/2}^{i+j}
V(k).
\end{equation}

\section{$N=1$ Neveu-Schwarz superalgebra and its minimal models}

The $N=1$ Neveu-Schwarz superalgebra is given by
$$\mathfrak{n}\mathfrak{s}=\bigoplus_{n\in \mathbb{Z}}\mathbb{C}L_{n}\bigoplus
\bigoplus_{n\in \mathbb{Z}}\mathbb{C}G_{n+1/2}\bigoplus \mathbb{C}C,$$
together with the following
$N=1$ Neveu-Schwarz relations:
\begin{eqnarray*} \label{nscom}
{[L_{m}, L_{n}]}&=&(m-n)L_{m+n}+\frac{C}{12}
(m^{3}-m)\delta_{m+n, 0},\\
{[L_{m}, G_{n+ 1/2}]}&=&\left(\frac{m}{2}-\left(n+
\frac{1}{2}\right)\right)G_{m+n+ 1/2},\\
{[G_{m+1/2}, G_{n-1/2}]}&=& 
2L_{m+n}+\frac{C}{3}(m^{2}+m)\delta_{m+n, 0},\\
{[C, L_{m}]}&=&0,\\
{[C, G_{m+1/2}]}&=&0
\end{eqnarray*}
for $m, n\in \mathbb{Z}$. 
We have the standard triangular decomposition
$\mathfrak{ns}=\mathfrak{ns}_+ \oplus \mathfrak{ns}_0 \oplus \mathfrak{ns}_-$
(cf. \cite{KWa}).
For every $(h,c) \in {\bf C}^2$, we denote  by
$M(c,h)$ Verma module for
$\mathfrak{ns}$ algebra.
For each $(p,q) \in {\bf N}^2$, $p = q \ {\rm mod}  \ 2$,
 let us introduce a family of complex 'curves'
$(h_{p,q}(t),c(t))$;
$$h_{p,q}(t)=\frac{1-p^2}{8}t^{-1}+\frac{1-pq}{4}+\frac{1-q^2}{8}t,$$
$$c(t)=\frac{15}{2}+3t^{-1}+3t.$$
Then from the determinant formula (see \cite{KWa} ) it follows
that $M(c,h)$ is reducible
if and only if there is a $t \in {\bf C}$ and $p,q \in {\bf N}$, $p= q  \ {\rm mod} \ 2$  such that
$c=c(t)$ and $h=h_{p,q}(t)$. In this case
$M(c,h)$ has a {\em singular vector} (i.e., a vector annihilated by
$\mathfrak{ns}_+$ ) of the weight $h+\frac{pq}{2}$.
Any such vector we denote by $v_{\frac{pq}{2}}$.

In this paper we are interested in the case $t=-1$. Then
$c(-1)=\frac{3}{2}$ and
$h_{p,q}(-1)=\frac{(p-q)^2}{8}$.
$h_{p,q}(-1)=h_{1,p-q+1}(-1)$, so we consider only
the case $h_{1,q}:=h_{1,q}(-1)$, (here $q$ is odd and positive).
Hence, each Verma module $M(\frac{3}{2},h_{1,q})$ is
reducible.

The following result easily follows from \cite{D} (or \cite{Aa})
and \cite{KWa}:
\begin{proposition}
For every odd $q$, $M(\frac{3}{2},h_{1,q})$ has the
following embedding structure
\begin{equation}
\ldots \rightarrow M\left(\frac{3}{2},h_{1,q+4}\right) \rightarrow
M\left(\frac{3}{2},h_{1,q+2}\right) \rightarrow
M\left(\frac{3}{2},h_{1,q}\right)
\rightarrow 0.
\end{equation}
Moreover, we have the following exact sequence:
\begin{equation}
0 \rightarrow M\left(\frac{3}{2},h_{1,q+2}\right) \rightarrow
M\left(\frac{3}{2},h_{1,q}\right) \rightarrow
L\left(\frac{3}{2},h_{1,q}\right) \rightarrow 0,
\end{equation}
where $L(\frac{3}{2},h_{1,q})$ is the corresponding
irreducible quotient.
\end{proposition}
\epf

Benoit and Saint-Aubin (cf. \cite{BSA}) found an explicit expression
for the singular vectors $P_{{\rm sing}}v_{1,q} \in M
\left(\frac{3}{2},h_{1,q}\right)$ that generates 
the maximal submodule:
\begin{equation}
\sum_{N;k_1,...,k_N} \sum_{\sigma \in S_N}  (-1)^{\frac{q-N}{2}}
c(k_{\sigma(1)},...,k_{\sigma(k)})
G(-k_1/2)\ldots G(-k_N/2)v_{1.q},
\end{equation}
where $S_N$ is a symmetric group on $N$ letters
and the first summation is over all the partitions of $q$ into the odd
integers $k_1,..,k_N$ and
$$c(k_{\sigma(1)},...,k_{\sigma(k)})=\prod_{i=1}^N {k_i-1 \choose
(k_i-1)/2 } \prod_{j=1}^{(N-1)/2} \frac{4}{\sigma_{2j} \rho_{2j}},$$
where $\sigma_j=\sum_{l=1}^j k_l$ and $\rho_j=\sum_{l=j}^N k_l$.

In the special case: $q=1$, $h_{1,1}=0$, 
$M(\frac{3}{2},0)$ has a singular vector
$G(-1/2)v$ which generate the maximal submodule.
By quotienting we obtain a {\em vacuum} module
$L(\frac{3}{2},0)=M(\frac{3}{2},0)/\langle G(-1/2)v_{3/2,0} \rangle$.

\section{$N=1$ superconformal vertex operator superalgebra and intertwining 
operators}

We use the definition of $N=1$ superconformal vertex operator superalgebra
(with and without odd variables) as 
in \cite{B} (cf. \cite{Kac}) 
and \cite{HM} (see also \cite{KW}).

Let $\varphi$ be a Grassman (odd) variable such that $\varphi^2=0$.
Every $N=1$ superconformal vertex operator superalgebra 
$(V, Y, \mathbf{1}, \tau)$ 
can be equipped with a structure of  $N=1$ superconformal vertex operator algebra
with an odd variable via
\begin{eqnarray*}
Y( \ ,(x,\varphi)) : V\otimes V&\to &V((x))[\varphi], \\
u\otimes v&\mapsto& Y(u, (x, \varphi)) v,
\end{eqnarray*}
where
$$Y(u, (x, \varphi))v=Y(u, x)v+\varphi Y(G(-1/2)u, x)v$$
for $u, v\in V$.

The same formula can be used in the case of modules
for the superconformal vertex operator superalgebra $(V, Y,
\mathbf{1}, \tau)$ (see \cite{HM}).

It is known (see \cite{KW}) that 
$V(c,0):=M(c,0)/ \langle G(-1/2)v_{c,0} \rangle$ 
\footnote{We write $L(c,0)$ if $V(c,0)$ is
irreducible.} is
a $N=1$ superconformal vertex operator superalgebra.
Also, every lowest weight $\mathfrak{ns}$--module with the central charge
$c$, is a $V(c,0)$--module.
If $c=\frac{3}{2}$ then $V(\frac{3}{2},0)=L(\frac{3}{2},0)$. Hence
\begin{proposition}
Every irreducible
$L(\frac{3}{2},0)$--module
is isomorphic to 
$L(\frac{3}{2},h)$, for some $h \in {\bf C}$.
\end{proposition}
{\em Proof:} The proof is essentially the same as the one in Proposition 3.1.
\epf

Among all irreducible $L(\frac{3}{2},0)$--modules
we distinguish modules isomorphic to $L(\frac{3}{2},h_{1,q})$, 
$ q \in 2{\bf N}-1$. These representations we call {\em degenerate 
minimal models}.

\subsection{Intertwining operators and its matrix coefficients}

The notation of an intertwining operators for $N=1$ superconformal vertex
operator algebras is introduced in \cite{KW} and \cite{HM}.

Let $W_{1}$, $W_{2}$ and $W_{3}$ be a triple of $V$--modules
and $\mathcal{Y}$ an intertwining operator of type ${W_3
\choose W_{1} \ W_{2}}$. 
Then we consider the corresponding intertwining operator   
with an odd variable (cf. \cite{HM}):
\begin{eqnarray*}
\mathcal{Y}( \ , (x,\varphi)) \ : W_{1}\otimes W_{2}&\to &W_{3}\{x\}[\varphi]\\
w_{(1)}\otimes w_{(2)}&\mapsto& \mathcal{Y}(w_{(1)}, (x, \varphi)) w_{(2)},
\end{eqnarray*}
such that 
$$\mathcal{Y}(w_{(1)}, (x, \varphi)) w_{(2)}
=\mathcal{Y}(w_{(1)}, x) w_{(2)}+
\varphi\mathcal{Y}(G(-1/2)w_{(1)}, x) w_{(2)}.$$

Let $w_1$ be a lowest weight vector for the Neveu-Schwarz algebra of
the weight $h$.
From the Jacobi identity we derive the following formulas:
\bea \label{evencom}
&& [L(-n),{\mathcal Y}(w_1,x_2)]=(x_2^{-n+1}\frac{\partial}{\partial
x_2}+(1-n)h){\mathcal Y}(w_1,x_2), \nn
&&  [G(-n-1/2),{\mathcal Y}(w_1,x_2)]=x_2^{-n}{\mathcal Y}(G(-1/2) w_1,x_2), \nn
&& [L(-n),{\mathcal Y}(G(-1/2)w_1,x_2)]=(x_2^{-n+1}\frac{\partial}{\partial
x_2}+(1-n)(h+\frac{1}{2}){\mathcal Y}(G(-1/2)w_1,x_2), \nn
&&  [G(-n-1/2),{\mathcal Y}(G(-1/2)w_1,x_2)]=
(x_2^{-n}\frac{\partial}{\partial x_2}-2nhx_2^{-n-1}){\mathcal
Y}(w_1,x_2).
\eea
In the odd formulation we obtain
\begin{eqnarray} \label{oddcom}
&& [L(-n),\mathcal{ Y}(w_{1},(x_2,\varphi))] \nn
&&= (x_2^{-n+1}\partial_{x_2}+(1-n)x_2^{-n}(h
+1/2\varphi \partial_{\varphi}))\mathcal{ Y}(w_{1},(x_2,\varphi)) \nn
&& [G(-n-1/2),\mathcal{Y}(w_{1},(x_2,\varphi))] \nn
&&= (x_2^{-n}(\partial_{\varphi}-\varphi \partial_{x_2})-2n
x_2^{-n-1}(h \varphi)\mathcal{ Y}(w_{1},(x_2,\varphi)),
\end{eqnarray}
where $\partial_{\varphi}$ is the odd (Grassmann) derivative.

\subsection{Even and odd intertwining operators}

In \cite{HM} we proved that every intertwining operator
$$\mathcal{Y} \in I {L(c,h_3) \choose L(c,h_1) \ L(c,h_2)}$$
is uniquely determined by the operators ${\mathcal Y}(w_1,x)$
and ${\mathcal Y}(G(-1/2)w_1,x)$, where $w_1$ is the lowest weight
vector of $L(c,h_1)$. This fact will be used later in connection 
with the following definition. 

\begin{definition} \label{oddeven}
{\em Let $| \ |$ denote the ($\mathbb{Z}/2 \mathbb{Z}$--valued) 
parity operator from the union of odd and even subspaces
for $V$--modules $W_i$, $i=1,2,3$.
An intertwining operator ${\mathcal Y} \in I {W_3 \choose W_1 \ W_2}$ is:
\begin{itemize}
\item {\em even}, if 
$$|{\rm Coeff}_{x^s} {\mathcal Y}(w_1,x)w_2|=|w_1|+|w_2|,$$
\item {\em odd}, if
$$|{\rm Coeff}_{x^s} {\mathcal Y}(w_1,x)w_2|=|w_1|+|w_2|+1,$$
\end{itemize}
for every $s \in \mathbb{C}$ and 
every $\mathbb{Z}/2\mathbb{Z}$--homogeneous vectors $w_1$ and $w_2$.}
\end{definition}

The space of even (odd) intertwining operators of the type
${W_3 \choose W_1 \ W_2}$ we denote
by $I \ {W_3 \choose W_1 \ W_2}_{\rm even}$ ($I \ {W_3 \choose W_1 \ W_2}_{\rm
odd}$).
In general we do not have  a decomposition of 
$I {W_3 \choose W_1 \ W_2}$ into the even and odd subspaces.

\subsection{Frenkel-Zhu's theorem for vertex operator superalgebras}
According to \cite{KW} (after \cite{Z}),
to every vertex operator superalgebra we can associate 
the Zhu's associative algebra $A(V)$. 
If $V=L(c,0)$, $A(L(c,0)) \cong {\bf C}[y]$.
where $y=[(L(-2)-L(-1)){\bf 1}]=[L(-2){\bf 1}]$ 
(because of the calculations that follow it is convenient to use
$y=[(L(-2)-L(-1)){\bf 1}]$).
Also to every $V$--module 
$W$  we associate a $A(V)$--bimodule $A(W)$ (cf. \cite{KW}).
In a special case $W=M(c,h)$, we have
$$A(M_{\mathfrak{ns}}(c,h))=M_{\mathfrak{ns}}(c,h)/O(M_{\mathfrak{ns}}(c,h)),$$
where 
\begin{eqnarray} \label{bimodules}
&& O(M_{\mathfrak{ns}}(c,h))=  \{L(-n-3)-2L(-n-2)+L(-1)v, \nn 
&& G(-n-1/2)-G(-n-3/2)v, n \geq 0, v \in M(c,h) \}.
\end{eqnarray}
It is not hard to see that, as $\mathbb{C}[y]$--bimodule, 
$$A(M(c,h)) \cong {\bf C}[x,y] \oplus {\bf C}[x,y]v,$$
where $v=[G(-1/2)v_{h}]$ and
$$y=[L(-2)-L(-1)], \ \ x=[L(-2)-2L(-1)+L(0)].$$

Let $W_{1}$, $W_{2}$ and $W_{3}$ be three 
$\mathbb{N}/2$--gradable irreducible $V$--modules
such that ${\rm Spec}L(0)|_{W_i} \in h_i + \mathbb{N}$, $i=1,2,3$
and $\mathcal{Y} \in I \ {W_3 \choose W_1 \ W_2}$. We define
$o(w_1):={\rm Coeff}_{x^{h_3-h_1-h_2}}\mathcal{Y}(w_1,x)$.
Because the fusion rules formula in \cite{FZ}
needs some modifications (cf. \cite{L1}) the
same modification is necessary for the main Theorem in \cite{KW}
(this can be done with a minor super--modifications along the lines of
\cite{L1}). Nevertheless (cf. \cite{KW}):
\begin{theorem} \label{superfz}
The mapping 
$$\pi : I {W_3 \choose W_{1} \ W_{2}} \rightarrow
{\rm Hom}_{A(V)}(A(W_{1})\otimes _{A(V)}W_{2}(0), W_{3}(0)),$$
such that
\begin{equation}
\pi(\mathcal{Y})(w_1 \otimes w_2)=o(w_1)w_2,
\end{equation}
is injective.
\end{theorem}

\section{Some Lie superalgebra homology}

In this section we recall
some basic definition from the homology theory
of infinite dimensional Lie superalgebras
which is in the scope of the monograph \cite{F} (in the cohomology
setting though).

Let $\mathcal L$ be an any (possibly infinite dimensional) $\mathbb
{Z}/ 2\mathbb{Z}$--graded Lie
superalgebra with the $\mathbb{Z}/2\mathbb{Z}$--decomposition:
${\mathcal L}={\mathcal L}_0 \oplus {\mathcal L}_1$. 
and let $M=M_0 \oplus M_1$ be any ${\bf Z}_2$--graded
${\mathcal L}$--module, such that the gradings are compatible.
Then, we form a chain complex $(C,d,M)$ (for details see \cite{F}),
$$0 \stackrel{d_0}{\leftarrow} C_0({\mathcal L},M) \stackrel{d_1}{\leftarrow}
C_1({\mathcal L},M) 
\stackrel{d}{\leftarrow }\ldots ,$$
where 
$$C_q({\mathcal L},M)=\bigoplus_{q_0+q_1=q} M \otimes \Lambda^{q_0}{\mathcal L}_0 \otimes
S^{q_1}{\mathcal L}_1,$$
$$C_q^p({\mathcal L},M)=\bigoplus_{\stackrel{q_0+q_1=q }{q_1+r = p \ {\rm
mod} 2 } } 
M_r \otimes \Lambda^{q_0}{\mathcal L}_0 \otimes S^{q_1}{\mathcal L}_1,$$
for $p=0,1$.
The mappings $d$ are super--differentials.
For $q \in {\bf N}$ and $p=0,1$, we define $q$--th homology
with coefficients in $M$ as:
\begin{equation} \label{z2}
H_q^p({\mathcal L}, M)=
{\rm Ker}(d_q (C_q^p({\mathcal L},M)))_p/(d_{r+1}(C_{q+1}^p({\mathcal
L},M)))_p.
\end{equation}
In a special case $q=0$, we have
$$H_0^0({\mathcal L}, M)=M_0/({\mathcal L}_0 M_0 + {\mathcal L}_1 M_1),$$
and
$$H_0^1({\mathcal L}, M)=M_1/({\mathcal L}_1 M_0 + {\mathcal L}_0 M_1).$$

We want to calculate $H_q({\mathcal L}_{s},L(\frac{3}{2},h_{1,q}))$.
for the Lie superalgebra
$${\mathcal L}_s = \bigoplus_{n \geq 0} {\mathcal L}_s(n),$$
where $ {\mathcal L}_s(n)$ is spanned by the vectors
$L(-n-3)-2L(-n-2)+L(-n-1)$ and $G(-n-1/2)-G(-n-3/2)$
, $n \in {\bf N}$.
From (\ref{bimodules}) we see (cf. \cite{HM}) 
that $H_0({\mathcal L}_{s},M(c,h))$ is a 
$\mathbb{C}[y]$--bimodule such that:
\begin{equation}
H_0({\mathcal L}_{s},M(c,h)) \cong A(M(c,h)) \cong {\bf C}[x,y] \oplus
{\bf C}[x,y]v.
\end{equation} 

\begin{remark}
{\em It is more involved to calculate $H_0(({\mathcal L}_{s},L(c,h))$, so we
consider only the special case $c=\frac{3}{2}$, $h=h_{1,q}$,
$q$ odd. As in the Virasoro case, it is easy to show that
the space $H_p({\mathcal L}_{s},L(\frac{3}{2},h_{1,q}))$ is infinite
dimensional for very $p,q,s \in {\bf N}$, and finitely generated as
a $A(L(3/2,0))$--module.
Moreover, it is not hard to see (by using the same method as in the Virasoro case) that
$${\rm Ext}^1_{ns, \mathcal{O}}(L(\frac{3}{2},h_{1,q}),L(\frac{3}{2},h_{1,r}))$$
is non--trivial (and one-dimensional) if and only if  $|r-q|=2$.}
\end{remark}

In the minimal models case we expect 
a substantially different result (cf \cite{FF1}). 
\begin{conjecture} \label{conj}
Let $c_{p,q}=\frac{3}{2}\left( 1-2\frac{(p-q)^2}{pq}\right)$ and
$h_{p,q}^{m,n}=\frac{(np-mq)^{2}-(p-q)^{2}}{8pq}$. Then
$${\rm dim} \ H_{q}({\mathcal L}_s, L(c_{p,q},h_{p,q}^{m,n}))< \infty,$$
for every $q \in {\bf N}$.
\end{conjecture}
There is strong evidence that Conjecture (\ref{conj})
holds based on \cite{Ad} and an example $c=-\frac{11}{14}$ treated in
Appendix of \cite{HM}.

The main difference between the
minimal models and the degenerate models is the fact that the maximal 
submodule for a minimal model is generated by two singular 
vectors, compared to 
$M(\frac{3}{2},h_{1,q})$ where the maximal submodule
is generated by a single singular vector.

\section{Benoit-Saint-Aubin's formula projection formulas} 
\subsection{Odd variable formulation}
We have seen before how to derive the commutation relation between
generators of $\mathfrak{ns}$ superalgebra and
${\mathcal Y}(w_1,x)$ where $w_1$ is a lowest weight vector for ${ns}$.
We fix ${\mathcal Y} \in I {L(\frac{3}{2},h) \choose
L(\frac{3}{2},h_{1,r}) \ L(\frac{3}{2},h_{1,q})}$ 
and consider the following matrix coefficient, 
\begin{equation}
\langle w'_{3}, \mathcal{Y}(w_1,x,\varphi)P_{{\rm sing}} w_2 \rangle,
\end{equation}
where $P_{sing}w_2=v_{1,q}$ (cf. (7),${\rm deg}(P_{sing})=q/2$) 
and $w_i$, $i=1,2,3$ are the lowest weight vectors. 

Since all modules are irreducible, by using
a result from \cite{HM} (Proposition 2.2), we get
$$\langle w'_{3}, \mathcal{Y}(w_1,x,\varphi)w_2 \rangle=c_1
x^{h-h_{1,q}-h_{1,r}} + c_2  \varphi x^{h-h_{1,q}-h_{1,r}-1/2},$$
where $c_1$ and $c_2$ are constants with the property 
\bea
&& c_1=c_2=0 \ {\rm implies} \ \ {\mathcal Y}=0.
\eea
From the formula (\ref{oddcom})
$$\langle w'_{3}, \mathcal{Y}(w_1,x,\varphi)P_{{\rm sing}} w_2 \rangle
=P(\partial_{x_2},\varphi) \langle w'_{3},
\mathcal{Y}(w_1,x,\varphi)w_2 \rangle,$$
where $P(\partial_{x_2},\varphi)$ is a certain
super-differential operator such that
$${\rm deg}(P_{{\rm sing}})={\rm deg}P(\partial_{x_2},\varphi)=q/2.$$
Therefore
$$P(\partial_{x_2},\varphi)c_1x^{h-h_{1,q}-h_{1,r}}=\varphi
C_1(h_{1,q},h_{1,r},h)x^{h-h_{1,q}-h_{1,r}-q/2},$$
and
$$P(\partial_{x_2},\varphi)\varphi c_2x^{h-h_{1,q}-h_{1,r}-q/2}=
C_2(h_{1,q},h_{1,r},h)x^{h-h_{1,q}-h_{1,r}-q/2}.$$

Constants $C_1(h_{1,q},h_{1,r},h)$ and 
$C_2(h_{1,q},h_{1,r}, h)$ (in slightly different form, but
in more general setting) were
derived in \cite{BSA}. 
Considering these
coefficients was motivated by deriving formulas for singular vectors from
already known singular vectors. 
By slightly modifying result from \cite{BSA} we obtain 
\begin{proposition} \label{main}
Suppose that ${\mathcal Y} \in I {L(\frac{3}{2},h) \choose
L(\frac{3}{2},h_{1,r}) \ L(\frac{3}{2},h_{1,q})}$ and
$P(\partial_x,\varphi)$ are as the above 
Then, up to a multiplicative constant,
$$C_1(h_{1,q},h_{1,r},h)=\prod_{-j \leq k \leq j} (h-h_{1,q+4k})$$
and
$$C_2(h_{1,q},h_{1,r},h)=\prod_{-j+1/2 \leq k \leq j-1/2}
(h+\frac{1}{2}-h_{1,q+4k}),$$
for $j=(r-1)/4$, $j >0$ (when $j=0$, $C_2(h_{1,1},h_{1,r},h)=1$).
\end{proposition}
{\em Proof:} 
The superdifferential operator $P(\partial_x,\varphi)$
is obtained by replacing generators $L(-m)$ and $G(-n-1/2)$
by the superdifferential operators
\begin{equation} \label{lm}
L(-m) \mapsto -(x_2^{-m+1}\partial_{x_2}+(1-m)x_2^{-m}(h_1
+1/2\varphi \partial_{\varphi}))
\end{equation}
and
\begin{equation}
G(-n-1/2) \mapsto (x_2^{-n}(\partial_{\varphi}-\varphi \partial_{x_2})-2n
x_2^{-n-1}(h_1 \varphi)),
\end{equation}
acting on $\langle w'_3, \mathcal{Y}(w_1,x,\varphi)w_2 \rangle$.
This action was calculated in \cite{BSA}. 
Their results (Formula 3.10 in
\cite{BSA}) implies the statement \footnote{In \cite{BSA}
a different sign was used in the equation (\ref{lm}). Still, we obtain
the same result if we consider an isomorphic algebra with the generators
$\tilde{L}(n):=-L(n)$. The same generators were used in \cite{FF2}.}.
\epf

\subsection{BSA formula without odd variables}
Since Frenkel-Zhu's formula does not involve odd variables
we need a version of Proposition \ref{main} without odd variables (which is of
course equivalent).
Again 
${\mathcal Y} \in  I \ {L(3/2,h) \choose L(3/2,h_{1,r}) \
L(3/2,h_{1,q})}$ is the same as the above. 
Then
$$\langle w'_{3}, \mathcal{Y}(w_1,x)P_{sing} w_2 \rangle=P_2(\partial_{x})
\langle w'_{3}, \mathcal{Y}(G(-1/2)w_1,x)w_2 \rangle,$$
and
$$\langle w'_{3}, \mathcal{Y}(G(-1/2)w_1,x)P_{sing} w_2 \rangle=P_1(\partial_{x})
\langle w'_{3}, \mathcal{Y}(w_1,x)w_2 \rangle,$$
where $P_1$ and $P_2$ are certain differential operators.
If 
$$P_2(\partial_{x})c_2x^{h-h_{1,q}-h_{1,r}-1/2}=c_2 K_2(h_{1,q}, h_{1,r},
h)x^{h-h_{1,q}-h_{1,r}-q/2 },$$
and 
$$P_1(\partial_{x})c_1x^{h-h_{1,q}-h_{1,r}}=c_1 K_1(h_{1,q}, h_{1,r}, h )
x^{h-h_{1,q}-h_{1,r}-q/2},$$
then, by comparing corresponding coefficients, we obtain 
\begin{eqnarray} \label{matching}
&& K_1(h_{1,q}, h_{1,r}, h)=C_1(h_{1,q}, h_{1,r}, h), \nn
&& K_2(h_{1,q}, h_{1,r}, h)=C_2(h_{1,q}, h_{1,r}, h).
\end{eqnarray}

Let us mention that the projection formulas from Proposition \ref{main}
have a simple explanation terms of {\em super density modules}
for the Neveu-Schwarz superalgebra.

\section{Fusion ring for the degenerate minimal models}

In order to obtain an upper bound for the fusion
coefficients (cf. Theorem \ref{superfz})
we first compute 
$$A(L(\frac{3}{2},h_{1,q})) \otimes_{A(L(3/2,0)}
L(\frac{3}{2},h_{1,r})(0).$$

${\mathbb Z}/2\mathbb{Z}$--grading of the $0$--th homology group (\ref{z2}) 
 enables
us (see Theorem (\ref{last}) to study odd and even intertwining
operators (see Definition \ref{oddeven}).
For that purpose we introduce the following splitting:
 
\begin{eqnarray} \label{zeroth}
&& A^0(L(\frac{3}{2},h_{1,q})):=H^0_0 ({\mathcal L}_s, L(\frac{3}{2},h_{1,q}))
\cong \frac{{\bf C}[x,y]}{I_1} \nn
&& A^1(L(\frac{3}{2},h_{1,q})):=H^1_0 ({\mathcal L}_s, L(\frac{3}{2},h_{1,q}))
\cong \frac{{\bf C}[x,y]v}{I_2},
\end{eqnarray}
where $I_1$ and $I_2$ are cyclic submodules (the maximal submodule
for $M(\frac{3}{2},h_{1,q})$ is cyclic !).
It seems hard to obtain explicitly these polynomials.
First we obtain some useful formulas
Inside $A(M(c,h))$ (cf. \cite{W}):
\bea \label{ok1}
&& [L(-n)v]=[((n-1)(L(-2)-L(-1))+L(-1))v]=\nn
&& [(n(L(-2)-L(-1))-(L(-2)-2L(-1)+L(0))+L(0))v]=\nn
&& (ny-x+{\rm wt}(v))[v].
\eea
for every $n \in \bf{N}$ and every homogeneous $v \in M(c,h)$.
Therefore in
$$A(M(\frac{3}{2},h_{1,q}))
\otimes_{A(L(\frac{3}{2},0))} L(\frac{3}{2},h_{1,r})(0)$$
we have
\begin{eqnarray} \label{ok32}
&& [L(-n)v]=(nh_{1,q}-x+L(0))[v]. \nn
&& [G(-n-1/2)v]=[G(-1/2)v].
\end{eqnarray}
Also, we have:
\bea \label{ok2}
&& [G(-n-\frac{1}{2})G(-m-\frac{1}{2})v]=[G(-1/2)G(-m-1/2)v]=\nn
&& [(2L(-m-1)-G(-m-1/2)G(-1/2))v]=[(2L(-m-1)-L(-1))v]=\nn
&& ((2m+1)y-x+{\rm wt}(v))[v].
\eea
By using (\ref{ok1}) and (\ref{ok2}) we obtain
\bea \label{ok3}
&& [G(-m_1-1/2) \ldots G(-m_{2r}-1/2)L(-n_1)...L(-n_s)v_{1,q}]=\nn
&& \prod_{i=1}^r ((2m_{2i}+1)h_{1,r}-x+\sum_{p=2i+1}^{2r}(m_p+1/2)+h_{1,q})
\cdot \nn
&& \prod_{j=1}^s (n_j h_{1,r}-x+\sum_{p=j+1}^s n_p+h_{1,q})[v].
\eea
inside
$$A(M(\frac{3}{2},h_{1,q}))
\otimes_{A(L(\frac{3}{2},0))} L(\frac{3}{2},h_{1,r})(0).$$
It is easy to obtain a similar formula for the vector
$$[G(-m_1-1/2) \ldots G(-m_{2r+1}-1/2)L(-n_1) \cdots L(-n_s)v_{1,q}].$$

\begin{lemma} \label{lelast}
Let $[P_{\rm sing}v_{1,q}]=Q_1(x)[G(-1/2)v_{1,q}]$ and
$[G(-1/2)P_{\rm sing}v_{1,q}]=Q_2(x)[v_{1,q}]$ be the projections inside
$$A(M(\frac{3}{2},h_{1,q}))
\otimes_{A(L(\frac{3}{2},0)} L(\frac{3}{2},h_{1,r})(0).$$
Then 
\bea
&& Q_1(h)=K_2(h_{1,q},h_{1,r},h), \nn
&& Q_2(h)=K_1(h_{1,q},h_{1,r},h),
\eea
for every $h \in \mathbb{C}$.
\end{lemma}
{\em Proof:} 
We use the notation from the section 6.2, where
$$\mathcal{Y} \in I \ {L(3/2,h) \choose L(3/2,h_{1,r}) \ L(3/2,h_{1,q})}.$$ 
By using (\ref{evencom}), we obtain
\bea \label{ok4}
&& \langle w'_3,
{\mathcal Y}(w_1,x)G(-m_1-1/2)...G(-m_{2r}-1/2)L(-n_1)...L(-n_s)w_2
\rangle= \nn
&& \prod_{i=1}^r -(x^{-m_{2i-1}-m_{2i}}\frac{\partial}{\partial x}
-2m_{2i}h_{1,r}x^{-m_{2i-1}-m_{2i}-1}) \cdot \nn
&& \prod_{j=1}^s
-(x^{-n_j+1}\frac{\partial}{\partial x}+(1-n_j)h_{1,r}x^{-n_j}) 
\langle w'_3, {\mathcal Y}(w_1,x)w_2 \rangle=\nn
&& c_1 \prod_{i=1}^r ((2m_{2i}+1)h_{1,r}-h+h_{1,q}+\sum_{p=2i+1}^{2r}
(m_{p}+1/2)) \cdot \nn
&& \prod_{j=1}^s (n_j h_{1,r}-h+\sum_{p=j+1}^s
n_p+h_{1,q})x^{h-h_{1,q}-
h_{1,r}-r-\sum m_i -\sum_j n_j},
\eea 
for the constant $c_1$ (see Section 6.1 and 6.2) that depends only on $\mathcal{Y}$.
There is a similar expression for
\begin{equation} \label{ok5}
\langle w'_3,
{\mathcal Y}(w_1,x)G(-m_1-1/2)...G(-m_{2r+1}-1/2)L(-n_1)...L(-n_s)w_2
\rangle.
\end{equation}
If we compare (\ref{ok3}) with (\ref{ok4}) (and corresponding
formulas for (\ref{ok5})) it follows that
$Q_1(h)$ is, up to a non--zero multiplicative constant, 
equal to $K_2(h_{1,r},h_{1,q},h)$ (singular vector is odd!)
and $Q_2(h)$ is, up to a multiplicative constant, equal to
$K_1(h_{1,r},h_{1,q},h)$.
\epf

Thus, Proposition \ref{main}
and Theorem \ref{last} gives us
\begin{theorem} \label{last}
\begin{itemize}
\item[(a)]
As a $A(L(3/2,0))$--module \\
\bea \label{amdec}
&& A(L(3/2,h_{1,q})) \otimes_{A(L(3/2,0))} L(3/2,h_{1,r})(0)
\cong \\
&& \frac{{\bf C}[x]}{<\prod_{-j \leq k \leq j}(x-h_{1,q+4k})>} \oplus 
\frac{{\bf C}[x]}{<\prod_{-j+1/2 \leq k \leq j+1/2} (h+1/2-h_{1,q+4k})>}. \nonumber
\eea
\item[(b)]
The space 
$$I \ {M(3/2,h)' \choose L(3/2,h_{1,q}) \ M(3/2,h_{1,r})},$$
is non--trivial if and only if $h=h_{1,s}$ for
some $s \in \{q+r-1,q+r-3,\ldots,q-r+1\}$.
\item[(c)]
The space 
$$I \ {L(3/2,h) \choose L(3/2,h_{1,q}) \ L(3/2,h_{1,r})},$$
is one--dimensional if and only if $h=h_{1,s}$, $s \in
\{q+r-1,q+r-3,\ldots,|q-r|+1\}.$
\end{itemize}
\end{theorem}
{\em Proof (a):}
From  Lemma \ref{lelast} it follows that
\bea
&& A(L(3/2,h_{1,r})) \otimes_{A(L(3/2,0))} L(3/2,h_{1,q}) \cong 
\frac{\mathbb{C}[x]}{\langle  Q_1(x) \rangle} \oplus \frac{\mathbb{C}[x]}{\langle  Q_2(x) \rangle}.
\eea
Now we apply (\ref{matching}) and Proposition \ref{main}.\\
{\em Proof (b):}
As in the Virasoro case, by examining carefully the main construction of 
intertwining operators in \cite{L1} with a minor super--modifications,
for every $A(L(3/2,0))$--morphism from 
$A(L(3/2,h_{1,q})) \otimes_{A(L(3/2,0))} L(3/2,h_{1,r})$ to $L(3/2,h)(0)$
we can construct a non--trivial intertwining operator of the form
$I \ {M(3/2,h)' \choose L(3/2,h_{1,q}) \ M(3/2,h_{1,r})}.$ \\
{\em Proof (c):} The proof and all the arguments involved are 
the same as in the Chapter 3. so we omit the
detials. We obtain a non--trivial intertwining operator
of the type
$ {L(3/2,h) \choose L(3/2,h_{1,q}) \ L(3/2,h_{1,r})}$
if $h=h_{1,s}$ for 
$s \in \{q+r-1,q+r-3,\ldots,q-r+1\} \cap
\{r+q-1,r+q-3,\ldots,r-q+1\}$, i.e.,
$s \in \{q+r-1,q+r-3,\ldots,|q-r|+1\}$.
\epf

\begin{theorem}
Suppose that $q \geq r$ \footnote{${W_3 \choose W_1 \ W_2} \cong 
{W_3 \choose W_2 \ W_1}$.}
\begin{equation}
{\rm dim} \ I \ {L(3/2,h_{1,s}) \choose L(3/2,h_{1,q}) \
L(3/2,h_{1,r})}_{\rm even}=1,
\end{equation}
if and only if
$$s \in \{ q+r-1,q+r-5,...,q-r+1 \}$$
\begin{equation}
{\rm dim} \ I \ {L(3/2,h_{1,s}) \choose L(3/2,h_{1,q}) \
L(3/2,h_{1,r})}_{\rm odd}=1
\end{equation}
if and only if
$$s \in \{ q+r-3,q+r-7,...,q-r+3 \}.$$
\end{theorem}
{\em Proof:}
By using (\ref{amdec}) we obtain 
the following decomposition:
\begin{eqnarray} \label{decomposition}
&& A^0(L(\frac{3}{2},h_{1,q})) \otimes_{A(L(3/2,0))} L(\frac{3}{2},h_{1,r})(0)
\cong \nn
&& {\mathbb C} v_{q+r-1} \oplus {\mathbb C} v_{q+r-5} \ldots \oplus {\mathbb C}v_{q-r+1}
\nn
&& A^1 L(\frac{3}{2},h_{1,q})) \otimes_{A(L(3/2,0))} L(\frac{3}{2},h_{1,r})(0)
\cong \nn
&& {\mathbb C} v_{q+r-3} \oplus {\mathbb C} v_{q+r-7} \ldots \oplus {\mathbb C}v_{q-r+3},
\end{eqnarray}
where ${\mathbb C}v_{i}$ is a ${\mathbb C}[y]$--module such that
$y.v_i=\frac{(i-1)^2}{8}v_i$. \\
{\em Claim:} Let 
$$\psi \in {\rm Hom}_{A(L(c,0))}(A^0(L(\frac{3}{2},h_{1,q}))
\otimes_{A(L(3/2,0))}L(\frac{3}{2},h_{1,r})(0),
L(\frac{3}{2},h_{1,s})(0)),$$
then the corresponding intertwining operator is even. 
Similarly if we start from
$$\psi \in {\rm Hom}_{A(L(c,0))}(A^1(L(\frac{3}{2},h_{1,q}))
\otimes_{A(L(3/2,0))}L(\frac{3}{2},h_{1,r})(0),
L(\frac{3}{2},h_{1,s})(0)),$$
the corresponding intertwining operator is odd.
\\
{\em Proof (of the Claim):} Let us elaborate the proof when $\psi$ is ``even''.
From the construction in \cite{FZ} and \cite{L2}
$\mathcal{Y}$ is obtained by lifting $\psi$ to a mapping from
$L(3/2,h_{1,q}) \otimes L(3/2,h_{1,r})(0)$ to
$L(3/2,h_{1,s})(0)$, such that
$$L(3/2,h_{1,q})_{\rm odd} \otimes L(3/2,h_{1,r})(0) \mapsto 0.$$
To extend this map to a mapping from
$L(3/2,h_{1,q}) \otimes M(3/2,h_{1,r})$ to $M(3/2,h_{1,s})'$
one uses generators and PBW so the sign is preserved. 
Because the isomorphism
$ I {W_3 \choose W_1 \ W_2} \cong I {W'_2 \choose W_1 W'_3}$
preserves the sign, i.e., odd intertwining operators
are mapped into odd and even into even, the result follows from the construction
of intertwining operators.
When $\psi$ is odd a similar argument works.
\epf

Let us summarize everything.
\begin{corollary} \label{last1}
Let ${\mathcal A}_s$ be a free abelian group with generators
$b(m), m \in 2{\bf N}+1$. 
Define a binary operation $\times: {\mathcal A}_s \times {\mathcal A}_s
\rightarrow {\mathcal A}_s$,
$$b(q) \times b(r)=\sum_{j \in {\bf N}} {\rm dim} \ I \ 
{L(3/2,h_{1,j}) \choose L(3/2,h_{1,q}) \
L(3/2,h_{1,r})} b(j).$$
Then ${\mathcal A}_s$ is a commutative associative
ring, and the mapping $b(m) \mapsto V(\frac{m-1}{4})$
gives an isomorphism to the representation ring
${\mathcal R}ep({\mathfrak osp}(1|2))$.
\end{corollary}
{\em Proof:} The proof follows from
Theorem \ref{last}(c) and (\ref{osp12}).
\epf

\section{Multiplicity $2$ fusion rules and super logarithmic intertwiners}

\subsection{A multiplicity $2$ case}
We have seen that in the $c=\frac{3}{2}$ case 
all fusion coefficients are $0$ or $1$.
Still, we expect (according to \cite{HM}) that for some
vertex operator superalgebras $L(c,0)$, fusion coefficients are $2$.

Here is one example.
If $c=0$, as in the case of the Virasoro algebra, the super vertex
operator algebra $L(0,0)=\frac{M(0,0)}{\langle G(-1/2)v_0,G(-3/2)v_0
\rangle}$ is trivial. 
Still we can consider a vertex operator superalgebra
$V(0,0):=\frac{M(0,0)}{\langle G(-1/2)v \rangle}$
Clearly, for every $h \in \mathbb{C}$,
we have (all modules are considered to be $V(0,0)$--modules):
\begin{equation}
{\rm dim} \ I \ {L(0,0) \choose L(0,h) \ L(0,h)}=2.
\end{equation}

The previous example is little bit awkward. 
Here is a nice example with irrational central charge:
\begin{proposition} \label{twons}
\begin{equation} \label{sqrt}
{\rm dim} \ I \
{L(\frac{15}{2}-3\sqrt{5},\frac{\sqrt{5}}{2}-1) \choose
L(\frac{15}{2}-3\sqrt{5},\frac{3}{4}(\frac{\sqrt{5}}{2}-1)) \
L(\frac{15}{2}-3\sqrt{5},\frac{3}{4}(\frac{\sqrt{5}}{2}-1))}=2.
\end{equation}
\end{proposition}
{\em Proof:}
It is not hard to see (by using a result form \cite{Aa} or \cite{D}) that 
$M(\frac{15}{2}-3\sqrt{5},\frac{3}{4}(\frac{\sqrt{5}}{2}-1))$
has the unique submodule that is irreducible (the case $II_+$ in
\cite{Aa}). 
If we analyze the determinant formula \cite{KWa}, singular vectors, 
and then use Theorem \ref{main}, we obtain
(\ref{sqrt}).
\epf

\subsection{A logarithmic intertwiner}
In \cite{M} we studied several examples of
logarithmic intertwining operators. Roughly, logarithmic intertwiners
exist if matrix coefficients yield some logarithmic solutions.
Our analysis can be extended for vertex operator superalgebras.
\begin{equation}
{\rm dim} \ I \ {W_2(\frac{27}{2},\frac{-3}{2}) \choose
L(\frac{27}{2},\frac{-3}{2}) \ L(\frac{27}{2},\frac{-3}{2})}=2,
\end{equation}
where $W_2(\frac{27}{2},\frac{-3}{2})$ is certain logarithmic module
(cf. \cite{M}).
The proof of this result and the discussion will appear in a separate
publication.

\section{Future work and open problems}
\begin{itemize}

\item
We know that  it is possible to
obtain {\em intertwining operator algebras} (see \cite{H2}) 
from the rational vertex operator algebras (satisfying
some natural convergence and extension condition and an
additional condition involving generalized modules).
Since the notation of intertwining operator
algebra can be (obviously) generalized such
that fusion algebra is an infinite--dimensional associative, commutative
algebra, one hopes that it is possible
to construct tensor categories for
degenerate minimal models.
In the language of conformal field theory this involves
explicit calculations of correlation functions for both 
products and iterates of intertwining operators (cf. Remark \ref{remaint}).

\item
Open problem: For rational vertex operator algebras, construct a {\em canonical}
isomorphism
$$A(M_1) \otimes_{A(V)} M_2(0) \cong
A(M_2) \otimes_{A(V)} M_1(0).$$

\item
(N=1 case) For which triples $L(c,h_1)$, $L(c,h_2)$ and $L(c,h_3)$ do
we have
$${\rm dim} \ I \ {L(c,h_3) \choose L(c,h_1) \ L(c,h_2)}=2 \ ?$$
\item
Determine the fusion ring for degenerate
minimal models for $N=2$ superconformal algebra by using our method
(it should be related to ${\rm Rep}(\mathfrak{osp}(2|2))$.

\item
Construct an analogue of the vertex tensor categories 
constructed in \cite{HM} (by using the main result in \cite{Ad}), for the
models studied in this paper.

\end{itemize}

\noindent {\small \sc Department of Mathematics, Rutgers University,
110 Frelinghuysen Rd., Piscataway, NJ 08854-8019}
\vskip 10mm 
\noindent {\small \sc Current Address: Department of Mathematics,
University of Arizona, Tucson, AZ 85721}

\noindent {\em E-mail address}: amilas@math.rutgers.edu, milas@math.arizona.edu


\begin{thebibliography}{FKRW}

\bibitem[AA]{Aa} 
A. Astashkevich, On the structure of Verma modules over
Virasoro and Neveu-Schwarz algebras, {\em Comm. Math. Phys.} 
{\bf 186}  (1997), 531--562.

\bibitem[A]{Ad}
D. Adamovi\'c, Rationality of Neveu-Schwarz vertex
operator superalgebra, {\em Internat. Math. Res. Notices} (1997), 865-874.

\bibitem[B]{B}
K. Barron, The supergeometric interpretation of vertex
operator superalgebras, Ph.D. thesis, Rutgers University, 1996.

\bibitem[BSA1]{BSA0}
L. Benoit and Y. Saint-Aubin,
Explicit expressions for some null vectors of the Virasoro algebra representations, 
{\em XVIIth International Colloquium on Group
Theoretical Methods in Physics}, Sainte-Adèle, PQ, 1988. 

\bibitem[BSA2]{BSA} 
L. Benoit and Y. Saint-Aubin,
Fusion and the Neveu-Schwarz singular vectors.
{\em Internat. J. Modern Phys.} A {\bf 9} (1994), 547--566. 

\bibitem[DG]{DG}
C. Dong, R. Griess,
Rank one lattice type vertex operator
algebras and their automorphism groups,  {\em 
J. of Algebra}, {\bf 208}, (1998), 262--275.

\bibitem[DL]{DL}
C. Dong, J. Lepowsky,
{\em Generalized vertex algebras and relative
vertex operators },  Progress in Mathematics Vol.112, 1993.

\bibitem[DM1]{DM1}
C. Dong, G. Mason, 
Quantum Galois theory for compact Lie groups. {\em J. of Algebra}, {\bf 214} (
5 92--102.

\bibitem[DM2]{DM2}
C. Dong, G. Mason
On quantum Galois theory. {\em Duke Math. J.}, {\bf 86}
(1997), 305--321.

\bibitem[DLM]{DLM}
C. Dong, H. Li, G. Mason
Compact automorphism groups of
vertex operator algebras, {\em Internat. Math. Res. Notices}, {\bf 18} (1996),

\bibitem[D]{D}
V. Dobrev,
Multiplet classification of the indecomposable highest weight modules over the Neveu--Schwarz and
Ramond superalgebras, {\em Lett. Math. Phys.} {\bf 11} (1986),
225--234; {\bf 13} (1987), 260. 


\bibitem[FF1]{FF1}
B. L. Feigin and D. B. Fuks,
 Cohomology of some nilpotent subalgebras of the
Virasoro and Kac-Moody Lie algebras,
{\em  J. Geom. Phys.} {\bf  5} (1988) , 209--235. 

\bibitem[FF2]{FF2}
B. L. Feigin and D. B. Fuks,  Representation of the Virasoro
algebra, in {\em Representations of Infinite--dimensional 
Lie groups and Lie algebras}, Gordon and Breach, 1989.

\bibitem[FF3]{FF3}
B. L. Feigin, D. B. Fuchs, Verma modules over the Virasoro algebra,
{\em  Lecture Notes in Math.}, {\bf 1060}, 230-245.

\bibitem[FM]{FM} 
B. Feigin and M. Malikov,
Modular functor and representation theory of
$\widehat{\rm sl}\sb 2$ at a rational level, 
In Operads: Proceedings of Renaissance Conferences
, {\em Contemporary  Math.} {\bf  202} 357--405. 

\bibitem[F]{F}
 D. B. Fuks, {\em Kogomologii beskonechnomernykh algebr Li}
(in Russian), ``Nauka'', Moscow, 1984

\bibitem[FHL]{FHL}
I.~B. Frenkel, Y.-Z. Huang and J.~Lepowsky,
On axiomatic approaches to vertex operator algebras and modules, 
{\em Memoirs Amer. Math. Soc.} {\bf 104}, 1993.


\bibitem[FKRW]{FKRW}
E. Frenkel, V. Kac, A. Radul, W. Wang, 
$W_{1+\infty}$ and 
$W(gl_N)$ with central charge $N$, 
{\em Comm. Math. Phys.} {\bf 170} (1995), 337-357

\bibitem[FLM]{FLM}
I.~B. Frenkel, J.~Lepowsky, and A.~Meurman,
{\em Vertex operator algebras and the Monster},
Pure and Appl. Math., {\bf 134}, Academic Press, New York, 1988.

\bibitem[FZ]{FZ} 
I. B. Frenkel and Y.  Zhu, Vertex operator
algebras associated to representations of affine and Virasoro
algebras, {\em Duke Math. J.} {\bf 66} (1992),  123--168.

\bibitem[H1]{H1}
Y.-Z. Huang, Virasoro vertex operator algebras,
(non--meromorphic) operator product expansion and the tensor product
theory, {\em J. Alg.} {\bf 182} (1996), 201--234.

\bibitem[H2]{H2}
Y.-Z. Huang, Genus-zero modular functors and intertwining operator algebras.
{\em Internat. J. Math.} {\bf 9} (1998), 845--863.

\bibitem[HL1]{HL1}
Y.-Z. Huang, J. Lepowsky, 
A theory of tensor products for module categories for a
vertex operator algebra  I, II,
{\em Selecta Math. (N.S.)} {\bf 1} (1995), 699--756, 757--786.


\bibitem[HL2]{HL2}
Y.-Z. Huang and J. Lepowsky, Tensor products of modules for a vertex
operator algebra and vertex tensor categories, in:
     {\em Lie Theory and Geometry,
in honor of Bertram Kostant,}
ed. R. Brylinski, J.-L. Brylinski, V. Guillemin, V. Kac,
Birkh\"{a}user, Boston, 1994, 349--383.

\bibitem[HM]{HM}
Y.-Z. Huang and A. Milas, 
Intertwining operator superalgebras 
and vertex tensor categories
for superconformal algebras, I
, math.QA/9909039, to appear in {\em Comm. Contem. Math.}.

\bibitem[KA]{K}
A. Kent,  Projections of Virasoro singular vectors,
{\em  Phys. Lett.} B {\bf 278} (1992), 443--448. 

\bibitem[KV]{Kac}
V. Kac, {\em Vertex algebras for beginners}, University Lectures
Series, Vol. 10, Providence, 1998.


\bibitem[KR]{KR}
V. Kac and A. Raina
Bombay lectures on highest weight representations of infinite-dimensional Lie algebras,
{\em Advanced Series in Mathematical Physics}, Vol 2, World
Scientific, NJ, 1987.




\bibitem[KWa]{KWa}
V. Kac and M. Wakimoto,
Unitarizable highest weight representations of the
Virasoro, Neveu-Schwarz and Ramond algebras,
{\em Conformal groups and related symmetries: physical
results and mathematical background, Lecture Notes in Phys.}  {\bf 261},
345--371.

\bibitem[KW]{KW} 
V. Kac and W. Wang,  Vertex operator 
superalgebras and their representations, 
in: {\em Mathematical aspects of conformal and topological field theories and 
quantum groups, Contemp. Math.} {\bf 175}, 161--191.

\bibitem[L1]{L1}
H. Li, Representation theory and a tensor product theory for vertex
operator algebras, PhD thesis, Rutgers University, 1994

\bibitem[L2]{L2}
H. Li, Determining fusion rules by $A(V)$-modules and bimodules. 
{\em J. of Algebra} {\bf 212} (1999), 515--556.

\bibitem[M]{M}
A. Milas, Weak modules and logarithmic intertwining operators, 
to appear in {\em Contemporary Mathematics}.


\bibitem[W]{W}
W. Wang, 
Rationality of Virasoro vertex operator algebras. 
{\em Internat. Math. Res. Notices} {\bf 7} (1993), 197--211. 

\bibitem[Z]{Z} 
Y. Zhu,  Modular invariance of characters of
vertex operator algebras, {\em J. Amer.  Math. Soc.}  {\bf 9}
(1996), 237--302.

\end{thebibliography}
\end{document}